\documentclass[12pt]{article}
\usepackage[matrix,arrow,curve]{xy}
\usepackage{amssymb}
\usepackage{supertabular}
\usepackage{a4}
\newcommand{\ddate}{\bf July 12 2002}
\newtheorem{dummy}{anything}[section] 

\newtheorem{Lemma}[dummy]{Lemma}
\newtheorem{Proposition}[dummy]{Proposition}

\newtheorem{Definition}[dummy]{Definition}
\newtheorem{Example}[dummy]{Example}
\newtheorem{Examples}[dummy]{Examples}
\newtheorem{Remark}[dummy]{Remark}
\newtheorem{Problem}[dummy]{Problem}
\newtheorem{ccote}[dummy]{}

\newcommand{\preu}{\noindent{\sc Proof: \ }}
\newcommand{\sk}[1]{\vskip #1 mm}
\newcommand{\eqref}[1]{(\ref{#1})}
\newcommand{\hfl}[2]{\smash{\mathop{\hbox to 1 truecm{\kern
3pt\rightarrowfill\kern 3pt}}%
\limits^{\scriptstyle#1}_{\scriptstyle#2}}}
\newcommand{\cqfd}{\unskip\kern 6pt\penalty 500%
\raise -2pt\hbox{\vrule\vbox to10pt{\hrule width
4pt\vfill\hrule}\vrule}\smallskip}
\newcommand{\bbf}{{\mathbb F}}
\newcommand{\bbr}{{\mathbb R}}

\newcommand{\bbc}{{\mathbb C}}
\newcommand{\bbz}{{\mathbb Z}}
\newcommand{\bbq}{{\mathbb Q}}
\newcommand{\bbn}{{\mathbb N}}

\newcommand{\cala}{{\mathcal A}}

\newcommand{\calc}{{\mathcal C}}
\newcommand{\cald}{{\mathcal D}}
\newcommand{\cale}{{\mathcal E}}

\newcommand{\calg}{{\mathcal G}}
\newcommand{\calh}{{\mathcal H}}
\newcommand{\cali}{{\mathcal I}}

\newcommand{\calk}{{\mathcal K}}

\newcommand{\caln}{{\mathcal N}}

\newcommand{\calp}{{\mathcal P}}

\newcommand{\calu}{{\mathcal U}}

\newcommand{\pcirc}{\kern .7pt {\scriptstyle \circ} \kern 1pt}

\newcommand{\mun}{{-1}}
\newcommand{\rmm}{\ensuremath{\bbr^m}}
\newcommand{\rmpos}{\ensuremath{(\bbr_{>0})^m}}
\newcommand{\rmcr}{\ensuremath{\bbr_{\scriptscriptstyle\nearrow}^m}}
\newcommand{\zmcr}{\ensuremath{\bbz_{\scriptscriptstyle\nearrow}^m}}
\newcommand{\strat}{{\rm Str}}
\newcommand{\strmcr}{{\rm Str}(\bbr_{\scriptscriptstyle\nearrow}^m)}
\newcommand{\strmpos}{{\rm Str}( (\bbr_{>0})^m)}
\newcommand{\fpp}{\ensuremath{\hookrightarrow}}
\newcommand{\chrmcr}{{\rm Ch}(\bbr_{\scriptscriptstyle\nearrow}^m)}
\newcommand{\chrmpos}{{\rm Ch}( (\bbr_{>0})^m)}
\newcommand{\ch}{{\rm Ch}}
\newcommand{\sol}{{\rm Sol}}
\newcommand{\nua}[2]{\caln^{#1}_{#2}}
\newcommand{\bnua}[2]{\bar{\caln}^{#1}_{#2}}
\newcommand{\calkk}{\dot\calk}

\newcommand{\su}[1]{\Sigma_{#1}}
\newcommand{\suo}[1]{\Sigma_{#1}^{or}}
\newcommand{\ns}{{\rm NS}}
\newcommand{\symm}{{\rm Sym}_m}
\newcounter{compol}
\newcommand{\ppv}[4]{\stepcounter{compol} {\small\thecompol} &%
$\langle #4\rangle$&$(#1)$}
\newcommand{\cvv}[4]{&#1&#2&#3\\[2.1mm]}
\newcommand{\pvi}[7]{\stepcounter{compol}{\small\thecompol} & %%
$\langle #7\rangle$ & $(#4)$  & #1 & #2 & #3}
\newcommand{\ccvi}[4]{  & #3 & #4 \\[1.3mm]}

\title{The space of clouds in an Euclidean space}
\author{Jean-Claude HAUSMANN and Eugenio RODRIGUEZ}
\date{\noindent\ddate}
\begin{document} \maketitle

\begin{abstract}\noindent
We study the space $\nua{m}{d}$ of clouds
in $\bbr^d$ (ordered sets of
$m$ points modulo the action of the group of
affine isometries).
We show that $\nua{m}{d}$ is a
smooth space, stratified over a certain
hyperplane arrangement in $\bbr^m$.
We give an algorithm to list all
the chambers and other strata
(this is independent of $d$).
With the help of a computer, we obtain
the list of all the chambers for $m\leq 9$
and all the strata when $m\leq 8$.
As the strata are the product of
a polygon spaces with a disk,
this gives a classification of $m$-gon spaces
for $m\leq 9$. When $d=2,3$,
$m=5,6,7$
and modulo reordering, we show that the chambers
(and so the different generic polygon spaces) are
distinguished by the ring structure of
their ${\rm mod}\,2$-cohomology.
\end{abstract}
\maketitle

\section{Introduction}\label{intro}

Let $E$ be an oriented finite-dimensional Euclidean space.
Let $\nua{m}{E}$ be the space
of ordered sets of $m$ points in $E$, modulo the group of
rigid motions of $E$; more precisely,
$$\nua{m}{E}:=G(E)\backslash E^m,$$
where the Lie group $G(E)$ is the semi-direct product of the
translation group of $E$ by $SO(E)$, the group
of linear orientation-preserving isometries of $E$, and
the group $G(E)$ acts diagonally on $E^m$.
We shall occasionally consider the space
$\bnua{m}{E}=\bar{G}(E)\backslash E^m$, where
$\bar{G}$ is the group of all affine isometries of $E^m$.
Observe that $\bnua{m}{E}$ is a subspace of
$\nua{m}{E'}$ when $E$ is a proper subspace of $E'$.
An element of $\nua{m}{E}$
will be called a {\it cloud} of $m$ points in $E$.
(The letter $\caln$ stands for ``nuage'', meaning
``cloud'' in French.) We abbreviate
$\nua{m}{\bbr^d}$ to $\nua{m}{d}$.
Observe that $\nua{m}{E}$ is canonically homeomorphic to
$\nua{m}{d}$, when $d=\dim E$.

The space $\nua{m}{d}$ plays a natural role in
celestial mechanics, at least for $d=2$ or $3$
(see, for instance, \cite{AC}).
Moreover, its importance was recognized especially in
statistical shape theory, a subject
which has developed rapidly during the
last two decades
(see \cite{Sm} and \cite{KBCL}) for a history).
There, the space $\nua{m}{d}$ is called the
{\it size-and-shape space} and is denoted by $S\Sigma^m_d$
\cite[\S\,11.2]{KBCL}.
This terminology and notation emphasizes
that $\nua{m}{d}$ is the cone,
with vertex $\nua{m}{0}$, over the {\it shape space} $\Sigma^m_d$,
defined as the quotient of $\nua{m}{d}-\nua{m}{0}$ by the
homotheties. A great amount is known about the homotopy type
of shape spaces. For instance, in \cite{KBCL},
Kendall, Barden, Carne and Li
show that $\Sigma^m_d$ admits cellular decompositions
leading to a complete computation of its homology groups.

In this paper, we present an alternative
decomposition of the space $\nua{m}{d}$. It is
based on polygon spaces,
a subject which has also encountered a rich development
during the last decade, in connection with Hamiltonian geometry.
This approach is completely different from that of
statistical shape theory and this paper is
essentially self-contained.

\iffalse
The cone structure of $\nua{m}{d}$ is useful to us because
of the role played by the clouds of points
having integral distances to their barycentre.
Ignorant anyway of the works
in statistical shape theory till the last phase of the redaction
of this paper, we do not use them and start from scratch.
\fi

First of all, the point set topology of $\nua{m}{d}$
is well behaved and $\nua{m}{d}$ is endowed with
a smooth structure. More precisely, the translations
act freely and properly on $E^m$
with quotient diffeomorphic to the vector subspace
$\calk(E^m)=\{(z_1,\dots,z_m)\in E^m\mid\sum z_i=0\}$.
Being therefore the quotient of $\calk(E^m)$ by the action
of the compact group $SO(E)$,
the space of clouds $\nua{m}{E}$ is locally compact
(in particular Hausdorff). Classical invariant theory
provides a proper topological embedding
$\varphi$ of $\nua{m}{d}$ into an Euclidean space $\bbr^N$
(see \ref{cloudssmooth}).
This embedding makes $\nua{m}{E}$ a {\it smooth space},
i.e. a topological space together with an algebra
$\calc^\infty(X)$ of {\it smooth functions} (with real values):
those functions which are locally the composition of
$\varphi$ with a $\calc^\infty$-function on $\bbr^N$.
One can prove that $f\in\calc^\infty(\nua{m}{E})$ if and only if
$f\pcirc\pi$ is smooth on $E^m$, where
$\pi:E^m\to\nua{m}{E}$ is the natural projection
(Proposition \ref{smooth-inv}). Any subspace of
of a smooth space naturally inherits a smooth structure and,
together with
{\it smooth maps} (see \ref{smoothspacesandmaps}),
smooth spaces
form a category whose equivalences are called
{\it diffeomorphisms}. Finally, we mention that the space
$\nua{m}{3}$ has the special feature that
the smooth maps admit a Poisson bracket
(see \ref{poisson}).

Our main tool for stratifying the space
$\nua{m}{E}$ is the map $\ell:\nua{m}{E}\to\bbr^m$
defined on $\rho=(\rho_1,\dots,\rho_m)\in E^m$ by
$$\ell(\rho):=\big(|\rho_1-b(\rho)|,\dots,|\rho_m-b(\rho)|\big),$$
where $b(\rho)=\frac{1}{m}\sum\rho_i$ is the barycentre of $\rho$.
This map $\ell$ is continuous and is smooth on
$\ell^\mun(\rmpos)$,
the open subset of points $\rho\in\nua{m}{E}$ such that
no $\rho_i$ is equal to $b(\rho)$.

We shall prove in Section \ref{PthA} that the
critical points of $\ell$ are the one-dimensional
clouds $\bnua{m}{1}\subset\nua{m}{E}$. The space
of critical values is then an arrangement of hyperplanes
in $\bbr^m$ that we shall describe now.
Let $\underline{m}:=\{1,2,\dots,m\}$
and denote by $\calp(\underline{m})$
the family of subsets of $\underline{m}$.
For $I\in\calp(\underline{m})$,
let $\calh_I$ be the hyperplane of $\bbr^m$ defined by
$$\calh_I:=\{(a_1,\dots,a_m)\in\bbr^m\mid
\sum_{i\in I}a_i=\sum_{i\notin I}a_i\}.$$
We call these hyperplanes {\it walls}. They
determine a stratification
$\calh(\rmm)$ of $\rmm$, i.e. a filtration
$$\{0\} = \calh^{(0)}(\rmm)\subset
\calh^{(1)}(\rmm)\subset\cdots\subset
\calh^{(m)}(\rmm)=\rmm,$$
with $\calh^{(k)}(\rmm)$ being
the subset of those $a\in\rmm$ which belong to at least
$m-k$ distinct walls $\calh_I$. A {\it stratum} of dimension
$k$ is a connected component of
$\calh^{(k)}(\rmm)-\calh^{(k-1)}(\rmm)$.
Note that a stratum of dimension $k\geq 1$ is an open convex cone
in a $k$-plane of $\rmm$.
Strata of dimension $m$ are called {\it chambers}.
We denote by $\strat(a)$ the stratum of $a\in\rmm$.
If $\strat(a)$ is a chamber, the $m$-tuple $a$
is called {\it generic} and we will often
denote $\strat(a)$ by ${\rm Ch}(a)$.

The stratification $\calh(\rmm)$ induces
a stratification of any subset $U$ of $\rmm$, in particular
of $U=\rmpos$.
Write $\strat(U)$ for the set of all the strata of $U$.
%One also gets, via the map $\ell$,
%a stratification $\calh(\nua{m}{E})$ by defining
%$\calh^{(k)}(\nua{m}{E}):=\ell^\mun(\calh^{(k)}(\rmm))$.
Denote by $\nua{m}{E}(a)$ the preimage
$\ell^\mun(\{a\})$ of $a\in\rmm$.
Strengthening results of
\cite{HK2}, we shall prove the following theorem in Section
\ref{PthA}.
\sk{2}\noindent {\bf Theorem A }\sl \
Let $a\in\rmpos$. Then there is a diffeomorphism
from $\ell^\mun(\strat(a))$ onto
to $\nua{m}{E}(a)\times\strat(a)$  intertwining the map $\ell$
with the projection to $\strat(a)$.
\rm \sk{2}
In the proof of Theorem A,
we actually construct,
when $\strat(a)=\strat(b)$, a diffeomorphism
$\psi_{ba}:\nua{m}{E}(a)\hfl{\approx}{}\nua{m}{E}(b)$.
One has $\psi_{ab}=\psi_{ba}^\mun$ and $\psi_{aa}={\rm id}$.
For $\alpha\in\strmpos$, we will sometimes use the notation
$\nua{m}{E}(\alpha)$ for any of the spaces $\nua{m}{E}(a)$
with $a\in\alpha$. This is in fact ambiguous because,
in general,
$\psi_{ca}\not=\psi_{cb}\pcirc\psi_{ba}$, so one cannot
use the maps $\psi_{ba}$ to define an equivalence relation
on $\ell^\mun(\strat(a))$ giving the points of
$\nua{m}{E}(\alpha)$. However,
$\psi_{ca}$ is isotopic to $\psi_{cb}\pcirc\psi_{ba}$,
and so the homotopy invariants of $\nua{m}{E}(\alpha)$,
for instance the elements of its cohomology ring,
are well defined.

Theorem A may provide good local models for
describing the evolution of a cloud. This is especially likely
when $\dim E=2,3$, where, for generic $a$, the
spaces $\nua{m}{E}(a)$ and thus $\ell^\mun(\strat(a))$
are smooth manifolds (see below).

Theorem A shows that $\nua{m}{E}$ is obtained
by gluing together pieces of the form
$\nua{m}{E}(\alpha)\times\alpha$
for various $\alpha\in\strmpos$.
From this point of view, the following questions are
natural.

\begin{enumerate}
\item Describe the set of all strata of $\rmpos$,
in particular the set of chambers.
This combinatorial problem does not depend on $E$.
\item Describe $\nua{m}{E}(\alpha)$ for all $\alpha\in\strmpos$.
\item Describe how a stratum of $\calh(\nua{m}{E})$
is attached to its bordering strata of lower dimension.
\end{enumerate}

The main issue of this paper is to answer Question 1 and, partly
Question 2 above. It is convenient to
take advantage of the right action of
the symmetric group $\symm$ on $\nua{m}{E}$ and on $\bbr^m$
by permutation of the coordinates
(to deal directly with the smooth space $\nua{m}{E}/\symm$
and get a corresponding statement of Theorem A, see
\ref{modsymsmooth}). This action permutes
the strata of $\calh(\rmpos)$, and $\nua{m}{E}(\alpha)$
is diffeomorphic to $\nua{m}{E}(\alpha^\sigma)$ for
$\sigma\in \symm$. The map $\ell$ is equivariant and
each $a\in\rmpos$ has a unique representative in $\rmcr$, where
$$\rmcr:=\{(a_1,\dots,a_m)\in\bbr^m\mid
0< a_1\leq\cdots\leq a_m\}.$$
Therefore, the set $\strmpos/{\rm Sym_m}$ is in bijection
with the set $\strmcr$.

In Sections \ref{Sgenchb} and \ref{Slowstr},
we show how to obtain a complete
list of the elements of $\chrmcr$ and $\strmcr$.
For this, we first show that the set of
inequalities defining a chamber $\alpha$ of $\rmcr$
can be recovered from
some very concentrated information which we call
the {\it genetic code} of $\alpha$. Abstracting some
properties of these genetic codes gives rise to
the combinatorial notion of a
{\it virtual genetic code}.
We design an algorithm to find all
virtual genetic codes, with the help of a computer
(the program in ${\rm C}^{++}$ is available at \cite{HRWeb}).
Deciding which virtual genetic code is the genetic code
of a chamber (realizability) is essentially done using
the simplex algorithm of linear programming.
We thus obtain the
list of all the chambers of $\rmcr$,
with the restriction $m\leq 9$ due to the computer's
limited capacities.
The set
$\strat(\bbr_{\scriptscriptstyle\nearrow}^{m-1})$
is determined using an injection of
$\strat(\bbr_{\scriptscriptstyle\nearrow}^{m-1})$
into $\chrmcr$ (see \S\,\ref{Slowstr}).
The number of elements of these sets is
\vskip 1mm
\stepcounter{dummy}
\begin{equation}\label{tablech}
\begin{minipage}{110mm}
\begin{tabular}{c|cccccccccccc} \small
$m$ &\footnotesize 3&\footnotesize 4&\footnotesize 5&\
\footnotesize 6 &\footnotesize 7&\footnotesize 8&\footnotesize 9
\\[1mm] \hline
\small $|\chrmcr|$ &\footnotesize 2 &\footnotesize 3&
\footnotesize 7&\footnotesize 21&\footnotesize 135&
\footnotesize 2470&  \footnotesize  175428&
\\[1mm] \hline
\small $|\strmcr|$
& \footnotesize 3& \footnotesize 7& \footnotesize 21
& \footnotesize 117& \footnotesize 1506
& \footnotesize 62254 & \footnotesize ? \\  %143696 not
\hline
\end{tabular}
\end{minipage}
\end{equation}\vskip 3mm

It turns out that, for $m\leq 8$, all virtual
genetic codes are realizable, but not for $m=9$:
only 175428 out of 319124 are realizable. The non-realizable
ones might well be of interest (see
Remark \ref{cohnonrea}).

Our algorithms produce, in each chamber $\alpha$,
a distinguished element $a_{\min}(\alpha)\in\rmcr$
with integral coordinates
and with $\sum a_i$ minimal. Several theoretical questions about
these elements $a_{\min}(\alpha)$ remain open (see \S\,\ref{Sgenchb}).

To describe the spaces $\nua{m}{E}(a)$ (Question 2 above),
we note that
$$\nua{m}{E}(a)=SO(E)\big\backslash
\{\rho\in E^m\mid \sum_{i=1}^m\rho_i=0 \hbox{ and } |\rho_i|=a_i\}.$$
The condition $\sum_{i=1}^m\rho_i=0$ suggests the picture of
a closed $m$-step piecewise-linear path in $E$,
whose $i$th step has length
$a_i$. Therefore, the space $\nua{m}{E}(a)$ is often called the
$m$-gon space (in $E$) of type $a$
(we could call it the space of clouds ``calibrated at $a$'').
These polygon spaces have been studied
in different notations, especially for $\dim E=2$ and 3 where,
for generic $a$, they are manifolds:
see, for instance, \cite{Kl}, \cite{KM}, \cite{HK1},
\cite{HK2}. For $\dim E>3$ or for $a$ non-generic,
see \cite{Ka1} and \cite{Ka2}.

The classification of
the polygon spaces $\nua{m}{E}(a)$,
for generic $a$, was previously known
when $\dim E=2,3$ and $m\leq 5$
(see, for instance \cite[\S\,6]{HK1}).
The genetic codes introduced in this paper
extend this classification up to $m=9$.
In \S\,\ref{Scalna}, we give handle-decomposition
information about the $6$-gon spaces $\bnua{6}{2}$
for the 21 chambers of $\bbr_{\scriptscriptstyle\nearrow}^6$.
This type of method could be applied to
any space $\nua{m}{E}(a)$ for generic $a$.
In addition to these geometric descriptions,
algorithms were previously found
which compute cohomological
invariants of the spaces $\nua{m}{3}(a)$,
for example their Poincar\'e polynomial
(\cite[Th. 2.2.4]{Kl}, \cite[Cor. 4.3]{HK2}).
This enables us, in Section \ref{Scoh},
to compute the Betti numbers of the spaces
$\nua{m}{3}(\alpha)$ for $m\leq 9$.
Moreover, presentations
of the cohomology ring of $\nua{m}{3}(\alpha)$
for any coefficients were given
in \cite[Th. 6.4]{HK2}.
This permits us to compute some invariants of
the ring $H^*(\nua{m}{3}(\alpha);\bbf_2)$
and prove in \ref{pfprB} the following result:

\sk{2}\noindent{\bf Proposition B }\sl\
For $5\leq m\leq 7$, the spaces $\nua{m}{3}(\alpha)$
(or $\bnua{m}{2}(\alpha)$),
for distinct chambers $\alpha$ of $\rmcr$,
have non-isomorphic ${\rm mod}\,2$-cohomology rings.
\rm\sk{2}

Here, the ring structure of $H^*(\nua{m}{3};\bbz_2)$
is important: the Betti numbers alone
do not distinguish the spaces.
Interestingly enough, the virtual genetic codes
which are not realizable also give rise
to non-trivial graded rings.
We do not know if these rings are cohomology rings of a space,
or of a manifold (see Remark \ref{cohnonrea}).

The paper is organized as follows. In Section \ref{Ssmooth},
we set the background of the smooth structure on $\nua{m}{E}$
which is used in Section \ref{PthA}
to prove Theorem A. In Section \ref{Sgenchb}, we introduce
the genetic code of a chamber and show how to obtain
the list of all chambers of $\rmcr$ for $m\leq 9$.
In Section \ref{Slowstr} we study the injection
$\strat(\bbr_{\scriptscriptstyle\nearrow}^{m-1})$
into $\strmcr$ and show how to obtain the list of
all strata of $\rmcr$ for $m\leq 8$. Section
\ref{Scalna} contains our information on
the spaces $\nua{m}{3}(a)$
and $\bnua{m}{2}(a)$ for generic $a$.
Section \ref{Scoh} is devoted the cohomology invariants
of the polygon spaces. Finally, the results of Sections
\ref{Scalna} and \ref{Scoh} are applied in Section \ref{hexasp}
to the case of hexagon spaces.

\paragraph{Acknowledgments:} Both authors thank the
Swiss National Fund for Scientific Research for its
support. We are indebted to R. Bacher for suggesting
the cohomology invariant $s(\alpha)$
of Proposition \ref{projvar}.

\section{The smooth structure on $\nua{m}{E}$}\label{Ssmooth}

\begin{ccote}\label{smoothspacesandmaps}
Smooth spaces and maps.\rm\
For $X$ a topological space, denote
by $\calc^0(X)$ the $\bbr$-algebra of continuous functions
on $X$ with real values.
If $h:X\to Y$ is a continuous map, denote by
$h^*:\calc^0(Y)\to\calc^0(X)$ the map
$h^*(f)=f\pcirc h$.

Let $X$ be a subspace of $\bbr^N$. A map $f:X\to\bbr$ is
{\it smooth} if, for each $x\in X$ there exists an open set
$U$ of $\bbr^N$ containing $x$ and a $\calc^\infty$
map $F:U\to\bbr$ which coincides with $f$ throughout $U\cap X$
(compare \cite[\S\,1]{Mi}). The smooth maps on $X$ constitute
a subalgebra $\calc^\infty(X)$ of $\calc^0(X)$.

More generally, if $\varphi:X\to\bbr^N$ is a topological embedding
of a space $X$ into $\bbr^N$ one may consider the subalgebra
$\calc^\infty(X)=\varphi^*(\calc^\infty(\varphi(X))$. We call
$\calc^\infty(X)$ a {\it smooth structure} on $X$ and $X$
(or rather the pair $(X,\calc^\infty(X))$) a
{\it smooth space}.

A continuous map $h:X\to Y$ between smooth spaces
is called {\it smooth}
if $h^*(\calc^\infty(Y))\subset\calc^\infty(X)$.
The map $h$ is a {\it diffeomorphism} if and only if it is a homeomorphism and
$h$ and $h^\mun$ are smooth. It is a {\it smooth embedding} if
$h^*(\calc^\infty(Y))=\calc^\infty(X)$. A smooth embedding is thus
a diffeomorphism onto its image.
\end{ccote}

\begin{ccote}\label{cloudssmooth}
The smooth structure on $\nua{m}{E}$.\rm\
Let $\kappa:E^m\to E^m$ be the linear projection
$$\kappa(z_1,\dots,z_m)=(z_1-b(z),\dots,z_m-b(z)),$$
where $b(z)=\frac{1}{m}\sum z_i$ is the barycentre of $z$.
The image of $\kappa$ is $\calk(E^m)$ and its
kernel is the diagonal $\Delta$ in $E^m$.

The normal subgroup $E$ in $G(E)$ of translations acts
freely and properly on $E^m$ and the quotient space
$E\backslash E^m$ is the same as the quotient vector space
$E/\Delta$. The projection $\kappa$ descends to a
linear isomorphism
$\bar\kappa:E\backslash E^m\hfl{\approx}{}\calk(E^m)$.
The space $\nua{m}{E}$ is now the quotient of $\calk(E^m)$
by the action of the compact group $SO(E)$.
Therefore $\nua{m}{E}$
is a locally compact Hausdorff space.

Consider the $m^2$ polynomial functions on $E^m$ given by
$z\mapsto\langle\kappa(z_i),\kappa(z_j)\rangle$,
where $\langle,\rangle$ denotes the scalar product
on $E$.
Choose an orientation on $E$. The determinants
$|\kappa(z_{i_1}),\cdots \kappa(z_{i_k})|$
with $i_1<\cdots<i_k$ ($k=\dim E)$ are another family of
$(^m_k)$ polynomial functions on $E^m$. All these
functions are $G(E)$-invariant and produce a continuous
map $\varphi:\nua{m}{E}\to\bbr^N$ with $N=m^2+(^m_k)$.
It is an exercise to prove that $\varphi$ is injective
and proper.
As $\nua{m}{E}$ is locally compact, the map $\varphi$
is a topological embedding of $\nua{m}{E}$ into $\bbr^N$
and its image is closed (for a family of inequalities defining
$\varphi(\nua{m}{E})$ as a semi-algebraic set, see \cite{PS}).

The embedding $\varphi$
endows $\nua{m}{E}$ with a smooth structure.
The following proposition identifies $\calc^\infty(\nua{m}{E})$
with the smooth functions on $E^m$ which are $G(E)$-invariant.
\end{ccote}

\begin{Proposition}\label{smooth-inv}
Let $f\in\calc^0(\nua{m}{E})$. The following are equivalent:
\begin{enumerate}
\renewcommand{\labelenumi}{(\Alph{enumi})}
\item $f\in\calc^\infty(\nua{m}{E})$.
\item There is a global $\calc^\infty$-function $F:\bbr^N\to\bbr$
such that $f=F\pcirc\varphi$.
\item  The map $f\pcirc\pi:E^m\to\bbr$ is $\calc^\infty$, where
$\pi:E^m\to\nua{m}{E}$ denotes the natural projection.
\end{enumerate}
\end{Proposition}

\preu
It is clear that (B) implies (A). Conversely,
let $f\in\calc^\infty(\nua{m}{E})$. For every $\rho\in\nua{m}{E}$,
one has an open set $U_\rho$ of $\bbr^N$ containing
$\varphi(\rho)$ and a smooth function $F_\rho:U_\rho\to\bbr$
with $f_\rho\pcirc\varphi=f$ on $\varphi^\mun(U_\rho)$.
Call $U_\infty=\bbr^N-\varphi(\nua{m}{E})$ and $F_\infty:U_\infty\to\bbr$
the constant map to $0$. As $\varphi(\nua{m}{E})$ is closed,
the family $\calu:=\{U_\rho\}_{\rho\in\nua{m}{E}\cup\{\infty\}}$
is an open covering of $\bbr^N$.
Let $\mu_\rho:\bbr^N\to\bbr$
be a smooth partition of the unity subordinated to $\,\calu$. Then
$F(x)=\sum_{\rho\in\nua{m}{E}\cup\{\infty\}}F_\rho(x)$ satisfies (B).

Statement (B) is obviously stronger than (C)
(which, incidentally, implies that $\pi$ is a smooth map).
For the converse,
one uses that the components of $\varphi$ constitute a
generating set for the algebra of
$SO(E)$-invariant polynomial functions on $\calk(E^m)$
\cite[\S\,II.9]{Wl}. Then,
any $SO(E)$-invariant smooth function
on $\calk(E^m)$ is of the form $F\pcirc\varphi$ by the Theorem
of G. Schwarz \cite{Sch}.
\cqfd

\begin{ccote}\label{cloudssmoothbar} \rm
The smooth structure on the space
$\bnua{m}{E}= G(E)\backslash E^m$ is obtained
as in \ref{cloudssmooth}. The embedding
$\varphi:\bnua{m}{E}\to\bbr^{m^2}$ is given by the
polynomial function $\rho\mapsto\langle\rho_i\rho_j\rangle$.
Proposition \ref{smooth-inv} holds true.
\end{ccote}

\begin{ccote}\label{modsymsmooth}
Clouds of unordered points.\ \rm
On $E^m=\{\rho:\underline{m}\to E\}$, the symmetric group
$\symm$ acts on the right, by pre-composition (or by permuting
the coordinates). This action descends on $\nua{m}{E}$.

As in \ref{cloudssmooth}, the space $\nua{m}{E}/\symm=G(E)\backslash E^m/\symm$
has a smooth structure, via a topological embedding
$\varphi:\bbr^m/\symm\to\bbr^N$ given by a generating
set of the algebra of polynomial functions
on $\calk(E^m)$ which are $SO(E)\times\symm$-invariant.
Proposition \ref{smooth-inv} holds true accordingly.

The space $\bbr^m/\symm$ also has a smooth structure via the
smooth embedding $\varphi:\bbr^m/\symm\to\bbr^m$
given by the $m$ elementary symmetric polynomials.
The map $\ell$ descends to
a continuous map $\bar\ell:\nua{m}{E}/\symm\to\bbr^m/\symm$
which is smooth away from $\ell^\mun(\{0\})$.
The composition $\psi:\rmcr\subset\rmpos\to\rmpos/\symm$ is
a smooth homeomorphism. The stratification
$\calh(\rmcr)$ can be transported via $\psi$ to
$\rmpos/\symm$, giving rise to a stratification
$\calh(\rmpos/\symm)$. The map $\bar\ell$ is stratified
and Theorem A holds true for $\bar\ell$.
Indeed, the diffeomorphisms constructed in
the proof of Theorem A given in \S\,\ref{PthA}
are natural with respect to the action of $\symm$.

We must be careful that the smooth homeomorphism
$\psi:\rmcr\to\rmpos/\symm$ is not
a diffeomorphism: the projection onto
the first coordinate is smooth on $\rmcr$ but not on
$\rmpos/\symm$.
\end{ccote}

\begin{ccote}\label{poisson}
Poisson structures on $\nua{m}{3}$.\ \rm
Recall that a Poisson structure on a smooth manifold $X$
is a Lie bracket
$\{,\}$ on $\calc^\infty(X)$ satisfying the Leibnitz rule:
$\{fg,h\}=f\{g,h\}+\{f,h\}g$. See \cite{MR}
for properties of Poisson manifolds. The same definition
makes sense on a smooth space.

The Euclidean space $E=\bbr^3$ has a standard smooth structure by
$$\{f,g\}(x):=\langle \nabla f \times \nabla g, x\rangle.$$
We endow the  product space $E^m$ with the product
Poisson structure. If $f,g:E^m\to\bbr$ are $SO(E)$-invariant,
so is the bracket $\{f,g\}$. Thus, the quotient space
$SO(E)\backslash E^m$ inherits a Poisson structure.

Using a canonical identification of $\bbr^3$ with $so(3)^*$,
the above Poisson bracket on $\bbr^3$
corresponds, up to sign, to the classical Poisson structure on $so(3)^*$
\cite[p. 287]{MR}.
The map $\mu: z\mapsto\sum_{i=1}^mz_i$, from $\bbr^3=so(3)^*$ to $E$
is the moment map for the diagonal action of $SO(E)$.
Let $\xi:\bbr^3\to\bbr$ be a linear map. By the Theorem of Noether
\cite[Th. 11.4.1]{MR}, if $f:E^m\to\bbr$ is a smooth
$SO(3)$-invariant map, then $\{f,\xi\pcirc\mu\}=0$.
This proves that $\{f,g\}=0$ for all $g\in\calc^\infty(E^m)$
such that $g_{|\calk(E^m)}=0$. Thus, the space $\nua{m}{3}$
inherits a Poisson structure so that the inclusion
$\nua{m}{3}\subset SO(E)\backslash E^m$
is a Poisson map.

When $a\in\rmpos$ is generic, the spaces $\nua{m}{3}(a)$
are manifolds and are the symplectic leaves
of $\ell^{\mun}(\strat(a))$.
This accounts for the symplectic structures on the polygon spaces
in $\bbr^3$ studied in \cite{Kl}, \cite{KM} and [HK1 and 2].
\end{ccote}

\section{Proof of Theorem A}\label{PthA}

Throughout this section,  the Euclidean space $E$ and the
number of points are constant.
Denote by $\calkk$ the subset of $m$-tuples
$\rho=(\rho_1,\dots,\rho_m)\in E^m$ such that $\rho_i\not=0$
and $\sum_{i=1}^m\rho_i=0$. Define the map
$\tilde\ell:\calkk\to\bbr^m$ by
$\tilde\ell(\rho):=(|\rho_1|,\dots,|\rho_m|)$.

An element $\rho=(\rho_1,\dots,\rho_m)\in E^m$ is called
{\it $1$-dimensional}
if the vector subspace of $E$ spanned by
$\rho_1,\dots,\rho_m$ is of dimension $1$
(therefore, $\rho$ represents an element of
$\bnua{m}{1}\subset\nua{m}{E}$).
These are precisely the singularities
of the map $\tilde\ell : \calkk\to \bbr^m$. Indeed:

\begin{Lemma}\label{ltha1}
Suppose that $\rho\in\calkk$ is not $1$-dimensional.
Then $T_\rho\tilde\ell$ is surjective.
\end{Lemma}

\preu
Let $(a_1,\dots,a_m)=\tilde\ell(\rho)$.
As $\rho$ is not
$1$-dimensional, there are two vectors
among $\rho_2,\dots,\rho_m$ which are
linearly independent.  The orthogonal
complements to these two vectors
then span $E$. Thus, there are curves $\rho_i(t)$ for $i=2,\dots,m$
such that $|\rho_i(t)|=a_i$ and
$$\sum_{i=2}^m\rho_i(t)=-(1+\frac{t}{a_1})\rho_1.$$
Therefore the map
$$t\mapsto ((1+\frac{t}{a_1})\rho_1,\rho_2(t),\dots,\rho_m(t))$$
represents a tangent vector $v\in T_\rho\calkk$ with
$T_\rho\tilde\ell(v)=(1,0,\dots,0)$.
The same can be done for the other
basis vectors of $\bbr^m$ proving that
$T_\rho\tilde\ell$ is surjective. \cqfd

Let $\rho\in\calkk$ be $1$-dimensional.
One thus has $\rho_i=\lambda_i\rho_m$ with $\lambda_i\in\bbr-\{0\}$.
Let $I(\rho)\in \calp(\underline{m})$ defined by
$i\in I(\rho)$ if and only if $\lambda_i<0$.
It is obvious that $\tilde\ell(\rho)$ belongs
to the wall $\calh_{I(\rho)}$.
\begin{Lemma} \label{ltha2}
Suppose that $\rho\in\calkk$ is $1$-dimensional.
Then the image of $T_\rho\tilde\ell$ is $\calh_{I(\rho)}$.
\end{Lemma}

\preu
Let $I=I(\rho)$.
One has $\sum_{i\in I}\rho_i=-\sum_{i\notin I}\rho_i$.
The components $\rho_i(t)$ of a curve
$\rho(t)\in E^m$ with $\rho(0)=\rho$
are of the form
$$\rho_i(t)=(1+\frac{c_i(t)}{a_i})\rho_i + w_i(t),$$
with $c_i(0)=0$ and $w_i(0)=0$,
where $c_i(t)\in\bbr$, and $w_i(t)$ is in the orthogonal complement
of $\rho_i$. The curve $\rho(t)$ is in $\calkk$ if and only if
$\sum_{i=1}^mw_i(t)=0$
and
\begin{equation}\label{ltha2eq1}
\sum_{i\in I}(1+\frac{c_i(t)}{a_i})\rho_i
=-\sum_{i\notin I}(1+\frac{c_i(t)}{a_i})\rho_i.
\end{equation}
Let $c(t)=(c_1(t),\dots,c_m(t))$.
The vector $\frac{\rho_i}{a_i}$ is constant when $i\in I$
and $\frac{\rho_i}{a_i}=-\frac{\rho_j}{a_j}$ if
$i\in I$ and $j\notin I$. Therefore,
Equation \eqref{ltha2eq1} is equivalent to $c(t)\in\calh_I$.
Finally, a direct computation shows that
the tangent vector $v\in T_\rho\calkk$ represented by $\rho(t)$
satisfies $T_\rho\tilde\ell(v)=\dot c(0)$. This proves the lemma.
\cqfd

\sk{2}
\noindent{\sc Proof of Theorem A : }
Let $a,b\in\rmpos$ be in the same stratum $\alpha$.
Let $X\subset\alpha$ be the segment joining $a$ to $b$.
For $\delta>0$, write $U_\delta:=\{x\in\bbr^m\mid d(x,X)<\delta\}$,
where $d(x,X)$ is the distance from $x$ to the segment $X$.
We choose $\delta$ small enough so that the walls meeting $U_\delta$,
if any, are only those containing $\alpha$.
Let $\tilde U_\delta:=\tilde\ell^\mun(U_\delta)\subset\calkk$.

Let $V^b$ be a vector field on $U_\delta$ of the form
$V^b_x=\lambda(x)(b-a)$, where $\lambda:U_\delta\to [0,1]$ is a
smooth function equal to $1$ on $U_{\delta/3}$ and to $0$ out of
$U_{2\delta/3}$.

Put on $\calkk$ and $\bbr^m$ the standard Riemannian metrics.
For $\rho\in\calkk$, define the vector subspace $\Delta_\rho$
of $T_\rho\calkk$ by
$\Delta_\rho:=(T_\rho\tilde\ell)^\sharp(T_{\ell(\rho)}(\bbr^m))$,
where $(T_\rho\tilde\ell)^\sharp$ is the adjoint of $T_\rho\tilde\ell$.
The vector spaces $\delta_\rho$ form a smooth distribution
(of non-constant rank) on $\calkk$.

The tangent map $T_\rho\tilde\ell$
sends $\Delta_\rho$ isomorphically onto the image of $T_\rho\tilde\ell$.
Since $X$ lies in $\alpha$,
Lemmas \ref{ltha1} and \ref{ltha2} show that $V^b_{\tilde\ell(z)}$
is in the image of $T_z\tilde\ell$ for all $z\in\tilde U_\delta$.
Therefore, there exists
a unique vector field $W^b$ on $\tilde U_\delta$ such that, for
each $z\in\tilde U_\delta$, one has $W^b_z\in\Delta_z$ and
$T_z\tilde\ell(W^b_z)=V^b_{\tilde\ell(z)}$.
The map $\tilde\ell$ being proper,
the vector field $W^b$ has compact support,
so its flow $\Phi_t$ is defined for all times $t$. Therefore,
$z\mapsto\Phi_1(z)$ gives a diffeomorphism
$\psi_{ba}:\tilde\ell^\mun(b)\hfl{\approx}{}\tilde\ell^\mun(a)$.

As its notation suggests, the map $\psi_{ba}$ depends only on $b$
and not on the choices involved in the definition of
$V^b$ ($\delta$ and $\lambda$). One can thus define
$\psi:\tilde\ell^\mun(\alpha)\to\nua{m}{E}(a)\times\alpha$
by $\psi(z):=(\tilde\ell(z),\psi_{\tilde\ell(z)a}(z))$.
The vector fields $V^b$ and $W^b$ depending smoothly on $b\in\alpha$,
the map is smooth as well as its inverse
$(x,u)\mapsto \psi_{ax}(u)$. Therefore, $\psi$ is a diffeomorphism.
As the Riemannian metric on $\calkk$ and the map $\tilde\ell$ are
invariant with respect to the action of $SO(E)\times\symm$
the map $\psi$ descends to a diffeomorphism
$\psi:\ell^\mun(\alpha)\hfl{\approx}{}\nua{m}{E}(a)\times\alpha$,
which proves Theorem A. Actually, each diffeomorphism
$\psi_{ba}$ descends to a diffeomorphism
$\psi_{ba}:\nua{m}{E}(b)\hfl{\approx}{}\nua{m}{E}(a)$.

\begin{Remark}\rm
Theorem A is also true for the spaces $\bnua{m}{E}$.
\end{Remark}

\section{The genetic code of a chamber}\label{Sgenchb}

Let $a\in (\bbr_{\geq 0})^m$.
Following \cite[\S\,2]{HK2}, we define
$S(a)\subset\calp(\underline{m})$ by
\stepcounter{dummy}\begin{equation}\label{shortdef}
I\in S(a) \ \Leftrightarrow\ \sum_{i\in I}a_i\leq\sum_{i\notin I}a_i.
\end{equation}
The very definition of the stratification $\calh$ implies that
$S(a)=S(a')$ if and only if $\strat(a)=\strat(a')$.
Thus, for $\alpha$ a stratum of $\rmpos$,
we shall write $S(\alpha)$ for the common set $S(a)$ with
$a\in\alpha$.

When $\alpha$ is a chamber, the inequalities in \eqref{shortdef}
are all strict. The elements of $S(\alpha)$ are then,
as in \cite[\S\,2]{HK2},
called  {\it short subsets} of $\underline{m}$. Observe that
$A\in\underline{m}$ is short if and only if its complement $\bar A$ is
not short. Therefore, if $\alpha$ is a chamber, the set $S(\alpha)$
contains $2^{m-1}$ elements.

Define $S_m(\alpha):=S(\alpha)\cap\calp_m(\underline{m})$,
where $\calp_m(\underline{m}):=
\{X\in \calp(\underline{m})\mid m\in X\}$.

\begin{Lemma}\label{Lcutmax}
Let $\alpha\in\chrmpos$.
Then $S(\alpha)$ is determined by $S_m(\alpha)$.
\end{Lemma}

\preu  One has
\begin{equation}\label{SmS}
I\in S(\alpha) \quad \Longleftrightarrow \quad
\left\{ \begin{array}{ccccc}
m\in I  \hbox{ and } I\in S_m(\alpha)\\[0mm]
\multicolumn{3}{l}{\hskip 15 mm\hbox{ or } }\\[0mm]
m\notin I  \hbox{ and }  \bar I\notin S_m(\alpha).  \cqfd
\end{array}\right.
\end{equation}
Let us now restrict ourselves to chambers of $\rmcr$.
We shall determine them by a very concentrated information
called their ``genetic code''.
Define a partial order ``\fpp'' on
$\calp(\underline{m})$ by saying that $A\fpp B$ if and only if there exits
a non-decreasing map $\varphi:A\to B$ such that $\varphi(x)\geq x$.
For instance $X\fpp Y$ if $X\subset Y$ since one can
take $\varphi$ being the inclusion.
The {\it genetic code} of $\alpha$ is the set of
elements $A_1,\dots,A_k$ of $S_m(\alpha)$ which are
maximal with respect to the order ``\fpp''.
By Lemma  \ref{Lcutmax}, the chamber $\alpha$ is determined
by its genetic code; we write
$\alpha=\langle A_1,\dots,A_k\rangle$ and call the sets $A_i$
the {\it genes} of $\alpha$.
Thanks to \eqref{SmS},
the explicit reconstruction of $S(\alpha)$ out of its
genetic code is given by the following recipe.

\begin{Lemma}\label{reconsalpha}
Let $\alpha=\langle A_1,\dots,A_k\rangle\in\chrmcr$.
Let $I\in\calp(\underline{m})$. Then
$$
I\in S(\alpha) \quad \Longleftrightarrow \quad
\left\{ \begin{array}{ccccc}
m\in I & \hbox{and} & \exists\, j\in\underline{k} \hbox{ with }
I\fpp A_j\\[1mm]
\multicolumn{3}{l}{\hskip 20 mm\hbox{ or } }\\[1mm]
m\notin I & \hbox{and} & \bar I\not\fpp A_j\ \forall\,
j\in\underline{k}. \cqfd
\end{array}\right.
$$
\end{Lemma}

\begin{Example}\label{3chbers}\rm
%Let us give some examples.
To unburden the notations, a subset $A$
of $\underline{m}$ is denoted by the number whose digits
are the elements of $A$ in decreasing order; example: $531=\{5,3,1\}$.
In $\bbr_{\scriptscriptstyle\nearrow}^3$, there are 2 chambers.
One of them, say $\alpha_0$,
contains points such as $(1,1,3)$ which are
not in the image of $\ell:\nua{3}{E}\to\bbr^3$.
One has $S_3(\alpha_0)=\emptyset$. Its genetic
code is empty and  one has
$$\alpha_0=\langle\rangle \quad ; \quad
S(\alpha_0) = \{\emptyset,1,2,21\}.$$
The other, $\alpha_1$ contains $(1,1,1)$, and one has
$$\alpha_1=\langle 3\rangle \quad ; \quad
S(\alpha_1) = \{\emptyset,1,2,3\}.$$
\end{Example}

Let us now figure out which subset
$\cala\subset\calp_m(\underline{m})$ is the genetic code of a
chamber of $\rmcr$. To reduce the number of trials,
observe that if $\alpha=\langle A_1,\dots,A_k\rangle$, then
\begin{enumerate}\renewcommand{\labelenumi}{(\alph{enumi})}
\item $A_i\not\fpp A_j$ for all $i\not=j$ and
\item $\bar A_i\not\fpp A_j$ for all $i,j$.
\end{enumerate}
Indeed, one has Condition (a) since the sets $A_i$ are maximal
(and we do not write them twice). For Condition
(b), if $\bar A_i\fpp A_j$, then $A_i$ would be
both short and not short and Inequalities \eqref{shortdef}
would have no solution. A finite set $\{A_1,\dots,A_k\}$,
with $A_i\in\calp_m(\underline{m})$ satisfying
Conditions (a) and (b) is called a {\it virtual genetic code}
(of type $m$), and we keep writing it by
$\langle A_1,\dots,A_k\rangle$.
Let $\calg_m$ be the set of
virtual genetic codes and $\calg_m^{(k)}$ the subset
of those  virtual genetic codes containing $k$ genes.

The determination of $\calg_m$ is algorithmic:
\begin{enumerate}
\item $\calg_m^{(0)}=\{\langle\rangle\}$.
\item Each $A\in\calp_m(\underline{m})$ satisfying
$\bar A\not\fpp A$ gives rise to a virtual genetic code
$\langle A\rangle$.
This gives the set $\calg_m^{(1)}$.
\item Suppose, by induction, that we know the set $\calg_m^{(k)}$
Then, each\\
$(\langle A_1,\dots,A_k\rangle,\langle A_{k+1}\rangle)$
in $\calg_m^{(k)}\times\calg_m^{(1)}$,
so that $\{A_1,\dots,A_{k+1}\}$ satisfies
Conditions (a) and (b),
gives rise to an element of $\calg_m^{(k+1)}$.
\end{enumerate}
When $\calg^{(k+1)}=\emptyset$, the process stops and
$\calg_m=\bigcup_{i=0}^m\calg_m^{(k)}$.

\begin{Examples}\rm
For $m=3$, the family $\calp_3(\underline{3})$ contains
the sets $3$, $31$, $32$ and $321$
(with the notations introduced in Example \ref{3chbers}).
Only $3$ satisfies
$\bar 3 = 21\not\fpp 3$. Thus $\calg_3^{(1)}=\{3\}$
while $\calg^{(2)}$ is empty.
We deduce that $\calg_3=\{\langle\rangle,\langle 3\rangle\}$.
They correspond to the two chambers of
$\bbr_{\scriptscriptstyle\nearrow}^3$ found in
Example \ref{3chbers}.
In the same way, we easily find the following
table:
\sk{2}
\begin{center}
\begin{tabular}{c|l}
$m$ & \hskip 10mm
Elements of $\calg_m$
\\ \hline
2 & $\langle\rangle$. \\[1pt]
3 & $\langle\rangle \, ,\, \langle 3\rangle$. \\[1pt]
4 & $\langle\rangle \, ,\, \langle 4\rangle \, ,\, \langle 41\rangle$.  \\[1pt]
5 & $\langle\rangle \, ,\, \langle 5\rangle \, ,\,
\langle 51\rangle \, ,\, \langle 52\rangle \, ,\,
\langle 53\rangle \, ,\, \langle 54\rangle \, ,\,
\langle 521\rangle$.\\
\hline\end{tabular}\end{center}
\end{Examples}

Having found the virtual genetic codes of type $m$, the next
question is which of them are
{\it realizable}, that is,
which of them is the genetic code of a chamber
of $\rmcr$. We proceed as follows. Each
virtual genetic codes $\langle A_1,\dots,A_k\rangle$
of type $m$ determines, a subset
$S_{\langle A_1,\dots,A_k\rangle}$ by the recipe
of Lemma \ref{reconsalpha}. Define the open polyhedral cone
$P:=P_{\langle A_1,\dots,A_k\rangle}$ by
$$P
:= \Big\{ x \in \rmcr \; \Bigm| \; \sum_{i \in I}x_i <
\sum_{i \notin I}x_i\; \forall \, I \in
 S_{\langle A_1,\dots,A_k\rangle} \Big\}.
$$
If there exits $\alpha\in\chrmcr$ with
$\alpha= \langle A_1,\dots,A_k\rangle$, then
$\alpha=P_{\langle A_1,\dots,A_k\rangle}$.
The realization problem is thus equivalent to
$P$ being not empty.
To find a point inside $P$,
we ``push'' its walls and consider:
\begin{equation}\label{defp1}
P_1:=
\Big\{ x \in \rmcr \; \Bigm| \;
\sum_{i \in I}x_i \leq \sum_{i \notin I}x_i
-1 \; \forall \, I \in S_{\langle A_1,\dots,A_k\rangle} \Big\} \subset
P
\end{equation}
As $P$ is an open cone in $\bbr^m$, then
$P$ is not empty if and only if $P_1$ is not empty.
Indeed, if $P$ is not empty, then
$\emptyset\not= P\cap\zmcr\subset P_1$.
We then use the {\it simplex algorithm} of linear programming
to minimize the $\ell_1$-norm $\sum_{i=1}^mx_i$ on $P_1$.
This algorithm either outputs an optimal solution,
which is a vertex of $P_1$, or
concludes that $P_1$ is empty \cite{Ch}.
%For details
%on the implementation of this algorithm see \cite{HRWeb}.

A program in $\rm C^{++}$ was designed,
following the above algorithms
(comments on this program and the source code can be
found in \cite{HRWeb}).
A computer could thus list all the chambers of $\rmcr$
for $m\leq 9$.
Each chamber $\alpha$ is given by a
distinguished element $a_{\min}(\alpha)\in\zmcr$ with
minimal $\sum_{i=1}^ma_i$.
The number of these chambers,
$|\chrmcr|=|\chrmpos/\symm|$
is the one given in the first line of
Table \eqref{tablech} in the introduction.

Experimentally, it turned out that, for $m\leq 8$, all virtual
genetic codes are realizable. This is not true for $m=9$:

\begin{Lemma}\label{notrel}
The virtual genetic code $\langle 9642\rangle\in\calg_9$
is not realizable.
\end{Lemma}

\preu
Let $S:=S_{\langle 9642\rangle}$.
As $9531\fpp 9642$, one has $9531\in S$. On the other hand,
$\overline{9642}=87531\notin S$. If $S=S(a)$ for some
generic $a\in\bbr_{\scriptscriptstyle\nearrow}^9$,
we would have
$a_7+a_8>a_9$. Now $965\notin S$ by Lemma \ref{reconsalpha},
therefore $\overline{965}=874321\in S$. By the above inequality
on the $a_i$'s, this would imply that
$94321\in S$ which contradicts $94321\not\fpp 9642$. \cqfd
%%Poincar\'e pol de $\langle 9642\rangle$:1,7,21,35,21,7,1

Our algorithm found 319124 elements in
$\calg_9$, out of which 175428 are realizable.

The list of all the chambers $\alpha$ of $\rmcr$
with their representative $a_{\min}(\alpha)$ can be found
further in this paper for $m\leq 6$
(Sections \ref{Scalna} and \ref{hexasp}) and on the
WEB page \cite{HRWeb} for $m=7,8,9$.

Several theoretical questions about $a_{\min}(\alpha)$ remain
open. For example, why $a_{\min}(\alpha)$ has integral
coordinates (with the $\ell_1$-norm
$|a|_1=\sum a_i$ odd)? {\it A priori}, the vertices of $P_1$
should only be in $\bbq_{\scriptscriptstyle\nearrow}^m$.
Is $a_{\min}(\alpha)$ always unique? This
suggests the following

\begin{ccote} {\sc Conjectures :}\sl

a) any stratum of $\alpha\in\calh(\rmcr)$ contains
a unique element $a_{\rm min}(\alpha)\in\zmcr$ with
minimal $\ell_1$-norm. \sk{1}

b) $\alpha$ is a chamber if and only $|a_{\rm min}(\alpha)|_1$
is an odd integer. \sk{1}

c) All vertices of $P_1(S)$ have integral coordinates.
\end{ccote}\rm

Conjecture b) is supported by the following evidences.
First, it is obvious that an element $a\in\zmcr$ with $|a|_1$ odd
is generic. On the other hand,
it is experimentally true for $m\leq 9$.
Conjecture a) for non generic strata is
experimentally true for $m\leq 8$ (see Section \ref{Slowstr}).
Conjecture c) has been checked for $m\leq 8$.

\begin{ccote}\label{cuts} Cuts:\ \rm
One can prove that the set $\calg_m$ of virtual genetic code
of type $m$ is in bijection with the set of ``{cuts}''
on $\underline{m}$ (the name is given in analogy
with the Dedeckind cuts of the rationals).
A subset $S$ of $\calp(\underline{m})$ is a {\it cut}
if, for all $I,J\subset \underline{m}$
the two following conditions are fulfilled:
\begin{enumerate}\renewcommand{\labelenumi}{(\Alph{enumi})}
\item  $I\in S \Leftrightarrow \bar I \notin S$.
\item  if $I\in S$ and $J\fpp I$, then $J\in S$.
\end{enumerate}
The bijection sends a cut $S$ of $\underline{m}$
to the set of maximal elements
(with respect to the order ``$\fpp$'') of $S_m$.
For details, see \cite{HRWeb}.
\end{ccote}

\section{Non generic strata}\label{Slowstr}

If $a\in\rmcr$ is not generic,
some inequalities of \eqref{shortdef} are equalities.
Thus, an element $I\in S(a)$ is either a {\it short subset}
of $\underline{m}$ (strict inequality) or an
{\it almost short subset}. As in Lemma \ref{Lcutmax},
$S(a)$ is determined by $S_m(a)=S(a)\cap\calp_m(\underline{m})$
and the latter is determined by those elements which are
maximal with respect to the order ``$\fpp$'' (the {\it genes}
of $S(a)$). We denote the genes
which are short subsets by $A_1,\dots,A_k$ and those
which are  almost short by $B_1^=,\dots,B_l^=$.
For instance, when $m=3$, one writes
$S(1,1,1)=\langle 3\rangle$ and $S(1,1,2)= \langle 3^=\rangle$.
To be more precise on our conventions, let $I^=$ be
an almost short gene of $S(a)$ and $J\fpp I$. If $|J|<|I|$, then $J$
is automatically short (since $a\in\rmcr\subset\rmpos$).
If $|J|=|I|$, then $J$ is supposed to be almost short unless there is
a short gene $K$ with $J\fpp K$. For instance,
$S(1,2,2,3,4)=\langle 51,53^=\rangle$.

The set $\strmcr$ of all the strata of
$\calh(\rmcr)$ will be studied via
a map $\alpha\mapsto\alpha^+$
from $\strat(\bbr_{\scriptscriptstyle\nearrow}^{m-1})$
to $\chrmcr$ which we define now.
Let $a\in\bbr_{\scriptscriptstyle\nearrow}^{m-1}$.
If $\varepsilon$ is small enough, the $m$-tuple
$a^+:=(\delta,a_1,\dots,a_{m-1})$ is a generic
element of $\rmcr$ for $\delta<\varepsilon$
and $\alpha^+:=\ch(a^+)$ depends
only on $\alpha=\strat(a)$.

If $\beta=\alpha^+$,
we denote $\alpha=\beta^-$. This makes sense because
of the following lemma.
\begin{Lemma}
The map $\alpha\mapsto\alpha^+$ is injective.
\end{Lemma}
\preu
Let $a,b\in\bbr_{\scriptscriptstyle\nearrow}^{m-1}$
such that $\strat(a)\not=\strat(b)$.
The segment joining $a$ to $b$ will then cross a wall $\calh_I$
which does not contain $\strat(a)$. But then, the
segment joining $a^+$ and $b^+$ will also
cross $\calh_I$, showing that $\ch(a^+)\not=\ch(b^+)$.
\cqfd

The correspondence $\alpha\mapsto\alpha^+$
can easily be described on the genetic codes.
The genetic code of $a^+$ has the same number of genes
than that of $a$. The correspondence goes as follows.
If $\{p_1\dots p_r\}$ is a gene of $S(a)$ which is short,
then $\{p_1^+\dots p_r^+,1\}$ is a gene of $S(a)^+$,
where $p_i^+=p_i+1$ (the genes of $S(a^+)$ are all short).
If $\{p_1\dots p_r\}^=$ is an almost short gene of $S(a)$,
then $\{p_1^+\dots p_r^+\}$ is a gene of $S(a^+)$.

The following convention will be useful.

\begin{ccote}\label{convzeros}\rm
Let $a=(a_1,\dots,a_k)$
be a generic element of $\bbz_{\scriptscriptstyle\nearrow}^k$.
For $m\geq k$, the $m$-tuple
$\hat a =(0,\dots,0,a_1,\dots,a_k)$ determines
a chamber $\hat\alpha$ represented by\\
$(\delta_1,\dots\delta_{m-k},a_1,\dots,a_k)$, where
$\delta_i>0$ and $\sum\delta_i<1$.
We say that  $\hat a$ is a {\it conventional representative}
of $\hat\alpha$. For instance,
$\langle 521\rangle$ having the conventional representative
$(0,0,1,1,1)$ shows that
$\langle 521\rangle=\langle 321\rangle^{++}$.
\end{ccote}

\begin{Lemma}
Let $\alpha$ be a chamber of $\rmcr$. Then
$\alpha=\beta^+$ if and only if one
(at least) of the two following statement holds:
\begin{enumerate}
\item $\alpha$ has a conventional representative $(0,a_2,\dots,a_m)$.
\item  there exists $(a_1,\dots,a_m)\in\alpha\cap\bbz^m$
with  $\sum a_i$ odd and $a_1=1$.
\end{enumerate}
\end{Lemma}

\preu
It is clear that either 1. or 2. implies $\alpha=\strat(a_2,\dots,a_m)^+$.
Also, if $\alpha=\beta^+$ for $\beta$ generic, then
$\alpha$ admits a conventional representative.
It remains to show that, if $\alpha=\beta^+$ with
$\beta$ non-generic, then Statement 2. holds true.

Observe that,
as the walls $\calh_I$ are
defined by linear equations with integral coefficients,
then $\beta\cap\bbq^{m-1}$ is dense in $\beta$.
As $\beta$ is a cone, it must contain
a point in $b\in\bbz_{\scriptscriptstyle\nearrow}^{m-1}$.
As $b$ is not generic, then
$\sum b_i$ must be even and the $m$-tuples $(\delta,b_1,\dots,b_{m-1})$
are all in $\alpha$ for $\delta<2$. \cqfd

Tables III-V of Section \ref{Scalna} and Table VI
of Section \ref{hexasp} show that $a_{\rm min}(\alpha)$
satisfies the above conditions for all
$\alpha\in\chrmcr$ when $m\leq 6$. This proves the following

\begin{Proposition}\label{plusbij}
The correspondence $\alpha\mapsto\alpha^+$ gives a bijection
$\strat(\bbr_{\scriptscriptstyle\nearrow}^{m-1})\to\chrmcr$
for $m\leq 6$.
\end{Proposition}

\iffalse
Proposition \ref{plusbij} together with
the tables mentioned above
give a classification of $\strmcr$ for $m\leq 5$.
\fi

Tables I and II below make the bijection $\alpha\mapsto\alpha^-$
from $\chrmcr$ to
$\strat(\bbr_{\scriptscriptstyle\nearrow}^{m-1})$ explicit (we put a conventional
$a_{min}(\alpha)$ when there exists one).

\vskip 4mm
\begin{center}
\begin{tabular}{ccccccc}
\multicolumn{6}{c}{\bf The bijection
$\ch(\bbr_{\scriptscriptstyle\nearrow}^4)\hfl{\approx}{}
\strat(\bbr_{\scriptscriptstyle\nearrow}^3)$}
\\ \hline\rule{0pt}{3ex}
$\alpha$ & $a_{\rm min}(\alpha)$ & $\nua{4}{2}(\alpha)$ & $\alpha^-$ &
$a_{\rm min}(\alpha^-)$ & $\nua{3}{2}(\alpha^-)$ \\[1pt]
\hline\rule{0pt}{3ex}
$\langle\rangle$ & $(0,0,0,1)$ & $\emptyset$ &
$\langle\rangle$ & $(0,0,1)$ & $\emptyset$ \\[2pt]
$\langle 4\rangle$ & $(1,1,1,2)$ & $S^1$ &
$\langle 3^=\rangle$ & $(1,1,2)$ & \footnotesize 1 point \\
$\langle 41\rangle$ & $(0,1,1,1)$ & $S^1{\scriptstyle\coprod}S^1$ &
$\langle 3\rangle$ & $(1,1,1)$ & \footnotesize 2 point \\
\hline
\end{tabular}
\end{center}

\vskip 4mm
\begin{center}
\begin{tabular}{ccccccc}
\multicolumn{6}{c}{\bf The bijection
$\ch(\bbr_{\scriptscriptstyle\nearrow}^5)\hfl{\approx}{}
\strat(\bbr_{\scriptscriptstyle\nearrow}^{4})$}
\\ \hline\rule{0pt}{3ex}
$\alpha$ & $a_{\rm min}(\alpha)$ & $\nua{5}{2}(\alpha)$ & $\alpha^-$ &
$a_{\rm min}(\alpha^-)$ & $\nua{4}{2}(\alpha^-)$ \\[1pt] \hline\rule{0pt}{3ex}
$\langle\rangle$ & $(0,0,0,0,1)$ & $\emptyset$ &
$\langle\rangle$ & $(0,0,0,,1)$ & $\emptyset$ \\[2pt]
$\langle 5\rangle$ & $(1,1,1,1,3)$ & $S^2$ &
$\langle 4^=\rangle$ & $(1,1,1,3)$ & \footnotesize 1 point\\[2pt]
$\langle 51\rangle$ & $(0,1,1,1,2)$ & $T^2$ &
$\langle 4\rangle$ & $(1,1,1,2)$ & $S^1$\\[2pt]
$\langle 52\rangle$ & $(1,1,2,2,3)$ & $\suo{2}$ &
$\langle 41^=\rangle$ & $(1,2,2,3)$ & $S^1\vee S^1$\\[2pt]
$\langle 521\rangle$ & $(0,0,1,1,1)$ & $T^2{\scriptstyle\coprod}T^2$ &
$\langle 41\rangle$ & $(0,1,1,1)$ & $S^1{\scriptstyle\coprod}S^1$\\[2pt]
$\langle 53\rangle$ & $(1,1,1,2,2)$ & $\suo{3}$ &
$\langle 42^=\rangle$ & $(1,1,2,2)$ &
\rule{0pt}{3ex}
\begin{minipage}{15 mm}
\setlength{\unitlength}{0.7mm}
\begin{picture}(0,0)(-10,0)
%\color{gristq}\graphpaper[2](0,-10)(180,100)\color{black}
\put(0,0){\circle{10}}
\put(-5,0){\circle*{1.5}}
\put(5,0){\circle*{1.5}}
\qbezier(-5,0)(0,5)(5,0)
\qbezier(-5,0)(0,-5)(5,0)
\end{picture}
\end{minipage}
\\[2pt]
$\langle 54\rangle$ & $(1,1,1,1,1)$ & $\suo{4}$ &
$\langle 43^=\rangle$ & $(1,1,1,1)$ &
\rule{0pt}{3ex}
\begin{minipage}{15 mm}
\setlength{\unitlength}{0.7mm}
\begin{picture}(0,0)(-10,0)
\put(0,0){\circle{10}}
\put(-2.5,4.3){\circle*{1.5}}
\put(-2.5,-4.3){\circle*{1.5}}
\put(5,0){\circle*{1.5}}
\qbezier(5,0)(0,0)(-2.5,4.3)
\qbezier(5,0)(0,0)(-2.5,-4.3)
\qbezier(-2.5,4.3)(0,0)(-2.5,-4.3)
\end{picture}
\end{minipage}
\\[5pt]
\hline
\end{tabular}
\end{center}\sk{3}

Here, $\suo{g}$ stands for the orientable surface of genus $g$
and the two graphs in the last column are the 2 and 3-fold
covers of $S^1\vee S^1$ without loops.
The same work with the bijection
$\ch(\bbr_{\scriptscriptstyle\nearrow}^6)\hfl{\approx}{}
\strat(\bbr_{\scriptscriptstyle\nearrow}^5)$ gives
the classification of all the 21 pentagon spaces (not necessarily generic)
obtained by A. Wenger \cite{We}.

In the above two tables, one sees that $\nua{5}{2}(\alpha)$
is the boundary of a
regular neighborhood (here in $\bbr^3$)
of $\nua{4}{2}(\alpha^-)$. This reflects the following fact.
Let $a_0\in\alpha$ without zero coordinate. For any
$a\in\alpha$, the Riemannian manifold $\nua{m}{E}(a)$
is canonically diffeomorphic to $\nua{m}{E}(a_0)$
by Theorem A and its proof. This produces
a family of Riemannian metrics $g_a$ on $\nua{m}{E}(a_0)$,
indexed by $a\in\alpha$.
When $a$ tends to a point $a^-\in\alpha^-$, the Riemannian manifold
$(\nua{m}{E}(a_0),g_a)$ converges, for the Gromov-Hausdorff metric, to the
metric space $\nua{m-1}{E}(a^-)$.

On the other hand, the map $\alpha\mapsto\alpha^+$ is not surjective
when $m\geq 7$. For instance, $\langle 764\rangle$,
with $a_{\rm min}=(2,2,2,2,3,3,3)$, is not of the
form $(\alpha^-)^+$. For, $\alpha^-$ would
be $\langle {653}^=\rangle$. As $421=\overline{653}\fpp 653$,
this would imply that all $a_i^-$ are equal
and $\alpha^- =\langle 654^=\rangle$. But
$\langle 654^=\rangle^+ = \langle 765\rangle \not= \langle 764\rangle$.

The table of $\ch(\bbr_{\scriptscriptstyle\nearrow}^7)$
(giving the 135 $7$-gon spaces) shows
18 chambers with the first coordinate $a_{min}$ not equal to $0$ or $1$
(see \cite{HRWeb}). One might ask whether there are other
$m$-tuples $a$ in these chambers with $a_1=0,1$. But, by
applying the simplex algorithm to minimize $a_1$ on the
polytope $P_1$ of \eqref{defp1}, we saw that this is not the case.
Therefore $|\strat(\bbr_{\scriptscriptstyle\nearrow}^6)|=118$.
The same procedure succeeded for $m=8$ and $9$, giving the
cardinality of $\strmpos/\symm=\strmcr$ for $m\leq 8$
listed in the introduction.

\section{Geometric descriptions of the $\nua{m}{2,3}(\alpha)$'s}\label{Scalna}

When $d=2$ or $3$ and $a$ is generic, the spaces
$\nua{m}{d}(a)$
are smooth manifolds,
since $SO(d)$ acts freely on the non-lined configurations.
The space
$\bnua{m}{2}(a)$ is also a manifold
and the map $\nua{m}{2}(a)\to\bnua{m}{2}(a)$
is a $2$-sheeted covering.
The space $\bnua{m}{2}(a)$
lies in $\nua{m}{3}(a)$ as
the fixed point set for the involution $\tau$ on
$\nua{m}{3}(a)$
obtained by reflection through a hyperplane.
Observe that $\dim\nua{m}{3}= 2(m-3)$ while
$\dim \dim\bnua{m}{2}=m-3$.
The manifold $\bnua{m}{2}$
plays the role of a real locus of
$\nua{m}{3}(a)$,
the latter being endowed with a natural Kaehler structure for
which the involution $\tau$ is antiholomorphic (see \cite[\S\,9]{HK2}).
It is shown in \cite[Th. 9.1]{HK2} that the cohomology rings
$H^{2*}(\nua{m}{3}(\alpha);\bbz_2)$
and $H^{*}(\bnua{m}{2}(\alpha);\bbz_2)$ are isomorphic,
by a graded ring isomorphism dividing the degrees by $2$.

The above polygon spaces were previously known for
$m\leq 5$ (see, for instance, \cite[\S\,6]{HK1}).
Our classification by genetic code
produces the more systematic tables below.
Conventional representatives $a_{\min}(\alpha)$
(see  \ref{convzeros}) are used when available.

\vskip 3mm
\begin{center}
\begin{tabular}{ccccccc}
\multicolumn{5}{c}{\bf Table III : the $3$-gon spaces}\\ \hline\rule{0pt}{3ex}
$\alpha$ & $a_{\rm min}(\alpha)$ & $\nua{3}{3}(\alpha)$
& $\bnua{3}{2}(\alpha)$ &
$\nua{3}{2}(a)$ \\[1pt] \hline\rule{0pt}{3ex}
$\langle\rangle$ & $(0,0,1)$ & $\emptyset$ & $\emptyset$ & $\emptyset$ \\[2pt]
$\langle 3\rangle$ & $(1,1,1)$ & \footnotesize 1 point &
\footnotesize 1 point &\footnotesize  2 points \\
\hline
\end{tabular}
\end{center}

\vskip 5mm
\begin{center}
\begin{tabular}{ccccccc}
\multicolumn{5}{c}{\bf Table IV : the $4$-gon spaces}\\ \hline\rule{0pt}{3ex}
$\alpha$ & $a_{\rm min}(\alpha)$ & $\nua{4}{3}(\alpha)$
& $\bnua{4}{2}(\alpha)$ &
$\nua{4}{2}(a)$ \\[1pt] \hline\rule{0pt}{3ex}
$\langle\rangle$ & $(0,0,0,1)$ & $\emptyset$ & $\emptyset$ & $\emptyset$ \\[2pt]
$\langle 4\rangle$ & $(1,1,1,2)$ & $\bbc P^1$ &
$\bbr P^1$ & $S^1$ \\[2pt]
$\langle 41\rangle$ & $(1,2,2,2)$ & $S^2$ &
$S^1$ & $S^1{\scriptstyle\coprod}S^1$ \\
\hline
\end{tabular}
\end{center}

\vskip 5mm
\begin{center}
\begin{tabular}{ccccccc}
\multicolumn{6}{c}{\bf Table V : the $5$-gon spaces}\\ \hline\rule{0pt}{3ex}
& $\alpha$ & $a_{\rm min}(\alpha)$ & $\nua{5}{3}(\alpha)$
& $\bnua{5}{2}(\alpha)$ &
$\nua{5}{2}(\alpha)$ \\[1pt] \hline
\rule{0pt}{3ex}
\ppv{0,0,0,0,1}{1}{1}{}\cvv{$\emptyset$}{$\emptyset$}{$\emptyset$}{}%x
\ppv{1,1,1,1,3}{7}{2}{5}\cvv{$\bbc P^2$}{$\bbr P^2$}{$S^2$}{}%x
\ppv{0,1,1,1,2}{5}{3}{51}%
\cvv{$\bbc P^2\,\sharp\,\overline{\bbc P}^2$}{$\su{1}$}{$T^2$}{}%x
\ppv{1,1,2,2,3}{9}{4}{52}%
\cvv{$(S^2\!\times\! S^2)\,\sharp\,\overline{\bbc P}^2$}{$\su{2}$}{$\suo{2}$}{}%x
\ppv{0,0,1,1,1}{3}{5}{521}\cvv{$S^2\times S^2$}{$T^2$}{$T^2{\scriptstyle\coprod}T^2$}{}%x
\ppv{1,1,1,2,2}{7}{6}{53}%
\cvv{$\bbc P^2\,\sharp\,3\overline{\bbc P}^2$}{$\su{3}$}{$\suo{3}$}{}%x
\ppv{1,1,1,1,1}{5}{7}{54}%
\cvv{$\bbc P^2\,\sharp\,4\overline{\bbc P}^2$}{$\su{4}$}{$\suo{4}$}{}%x
\hline
\end{tabular}
\end{center}

Our method produces a classification of the spaces
$\nua{m}{\bbr^n}(\alpha)$ for $m\leq 9$,
$\alpha$ a chamber, and $n\geq 2$.
Table VI of Section \ref{hexasp} gives the list of hexagon
spaces. The tables for generic $m$-gon spaces when $m=7,8,9$ are
too big to be included in this paper. They can be consulted
on the WEB page \cite{HRWeb}.

We shall now give procedures describing $\nua{m}{E}(\beta^+)$
in terms of $\nua{m-1}{E}(\beta)$ when $\beta$ is generic
and $\dim E=2,3$.
A $m$-tuple $(\rho_1,\dots,\rho_m)\in \calk(E^m)$ is called a
{\it vertical configuration} if $\rho_m=(0,\dots,0,-|\rho_m|)$.

\begin{Proposition}\label{timesS1} If
$\beta\in {\rm Ch}(\bbr_{\scriptscriptstyle\nearrow}^{m-1})$, then
$\nua{m}{2}(\beta^+)$ is diffeomorphic to
$\nua{m-1}{2}(\beta)\times S^1$.
\end{Proposition}

\preu  Let $(b_2,\dots,b_m)\in\beta$ and let $\varepsilon>0$
small enough so that $a:=(\varepsilon,b_2,\dots,b_m)\in\beta^+$.
A class in $\nua{m}{2}(a)$ has a unique representative
$\rho=(\rho_1,\dots,\rho_m)$ which is a vertical configuration.
As $b$ is generic, if $\varepsilon$ is
small enough, then $(b_2,\dots,b'_m)\in\beta$ when
$|b'_m-b_m|<\varepsilon$. The $(m-1)$-tuple
$(\rho_2,\dots,\rho_m+\rho_1)$ thus represents
an element $\rho'\in \nua{m-1}{2}(\beta)$ and
the correspondence $\rho\mapsto (\rho',\rho_1)$
produces a diffeomorphism from  $\nua{m}{2}(\beta^+)$ to
$\nua{m-1}{2}(\beta)\times S^1$. \cqfd

Let $\rho=(\rho_1,\dots,\rho_m)\in(\bbr^3)^m$.
Let $\rho_m^\perp$ be the orthogonal complement of $\rho_m$,
oriented by the vector $\rho$. Let $\pi:\bbr^3\to\bbr^2$ be the composition
of the orthogonal projection $\bbr^3\to\rho_m^\perp$ with some chosen
isometry $\rho_m^\perp\hfl{\approx}{}\bbr^2$ preserving the
orientation. If $a\in\rmpos$ is generic, than
$\pi(\rho_1)$, $\pi(\rho_1+\rho_2)$,\dots ,
$\pi(\rho_1+\cdots+\rho_{m-2})$
are not all zero. This defines
a smooth map
\begin{equation}\label{mapr}
r:\nua{m}{3}(\alpha)\to(\bbr^2)^{m-2}\!-\!\{0\}\big/SO(2)
\end{equation}
where $\alpha=\ch(a)$. The right hand member of Equation
\eqref{mapr} is homotopy equivalent to $\bbc P^{m-3}$. The
map $r$ thus determines a cohomology class
$R\in H^2(\nua{m}{3}(\alpha);\bbz)$ which is the characteristic
class of some principal circle bundle
$\cale(\alpha)\to\nua{m}{3}(\alpha)$.
The class $R$ was introduced in \cite[\S\,6 and 7]{HK2} and
will appear again in Section \ref{Scalna} below.

\begin{Lemma}\label{Ltotspace}
(compare \cite[Prop. 7.3]{HK2})
The total space $\cale(\alpha)$
is $S^1$-equivariantly diffeomorphic to the spaces
of representatives of $\nua{m}{3}(a)$,
($\ch(a)=\alpha$) which are vertical configurations.
\end{Lemma}

\preu  Let $\cale'(\alpha)\subset (\bbr^3)^m$ be the space described
in the statement.
Any element of $\nua{m}{3}(a)$ has at least one representative
which is a vertical configuration and any two of those are
in the same orbit under the orthogonal action of $S^1=SO(2)$
fixing the vertical axis.
As $a$ is generic, the quotient map
$\cale'(\alpha)\to\nua{m}{3}(a)$ is then a principal circle
bundle. If $\pi:\bbr^3\to\bbr^2$ denotes the projection onto
the first two coordinates, the correspondence
$\rho\mapsto \big(\pi(\rho_1),\pi(\rho_1+\rho_2),\dots ,
\pi(\rho_1+\cdots+\rho_{m-2})\big)$
defines a smooth $S^1$-equivariant
map $\tilde r:\cale'(\alpha)\to(\bbr^2)^{m-2}$
which covers the map~$r$. This proves that the characteristic
class of $\cale'(\alpha)\to\nua{m}{3}(a)$ is $R$. \cqfd

\begin{Example}\label{excut<m>}
\rm The chamber $\alpha=\langle m\rangle$ of $\rmcr$ has
its minimal representative
$a=a_{\min}(\alpha)=(1,\dots,1,m-2)$. As, in
a vertical configuration $\rho$ of
$[\rho]\in\nua{m}{3}(a)$,
one has $\sum_{i=1}^{m-1}\rho_i=(0,\dots,0,m-2)$,
the sequence of the third coordinate
of $\rho_1$, $\rho_1+\rho_2$, \dots, must be strictly increasing.
This implies that the map $\tilde r$ of
the proof of \ref{Ltotspace} is a smooth
$S^1$-equivariant embedding.
It induces
a diffeomorphism from $\nua{m}{3}(\langle m\rangle)$
onto $\bbc P^{m-3}$ and an identification of the
bundle $\cale(\langle m\rangle)\to\nua{m}{3}(\langle m\rangle)$
with the Hopf bundle.
(see also \cite[Remark 4.2]{Ha}).
\end{Example}

Let $\cald(\alpha)$ be the total space of the $D^2$-bundle
associated to  $\cale(\alpha)\to\nua{m}{3}(\alpha)$.

\begin{Proposition}\label{Pdouble}
 If $\beta\in {\rm Ch}(\bbr_{\scriptscriptstyle\nearrow}^{m-1})$,
then $\nua{m}{3}(\beta^+)$ is diffeomorphic to
the double of $\cald(\beta)$.
\end{Proposition}

\preu
Let $(b_2,\dots,b_m)\in\beta$ and let $\varepsilon>0$ be
small enough so that $a:=(\varepsilon,b_2,\dots,b_m)\in\beta^+$.
A class in $\nua{m}{3}(a)$ has a representative
$\rho=(\rho_1,\dots,\rho_m)$ which is a vertical configuration
and with $\rho_1=(\varepsilon\cos\theta,0,\varepsilon\sin\theta)$.
Let $\check\cale(a)$ be the space of such representatives.
%called {\it special vertical configurations}.
The map sending $\rho$ to $\theta$ is a smooth map
$\theta:\check\cale(a)\to [0,\pi]$.

If $\rho\in\check\cale(a)$, then
$\rho_m+\rho_1$ is close to $\rho_m$. This defines a
smooth map $P:\check\cale(a)\to SO(3)$, sending $\rho$
to $P_\rho$, characterized by
$P_\rho(\rho_m+\rho_1)=(0,0,-|\rho_m+\rho_1|)$ and
$P_\rho={\rm id}$ if $\theta(\rho)=0,\pi$.
The smooth map
$\check F:\check\cale(a)\to [0,\pi]\times \cale(\beta)$
given by
$$
\check F(\rho):=\Big(\,\theta(\rho)\, ,\,
(P_\rho(\rho_2),\dots,P_\rho(\rho_{m-1}),P_\rho(\rho_m+\rho_1))\,\Big)
$$
is a diffeomorphism. It induces a diffeomorphism
\begin{equation}\label{lastdiffeo}
F: \nua{m}{3}(\alpha)\cong\nua{m}{3}(a)\hfl{}{}
[0,\pi]\times \cale(\beta) \bigg/\sim
\end{equation}
where $\sim$ is the equivalence relation given
by $(0,\eta)\sim (0,g\cdot\eta)$
and $(\pi,\eta)\sim (\pi,g\cdot\eta)$
for all $g\in SO(2)$.
The right member of \eqref{lastdiffeo}
is diffeomorphic to the double of
$\cald(\beta)$ which proves the proposition. \cqfd

If, in \eqref{mapr}, one replaces $\bbr^3$ by $\bbr^2$, one gets
a map
$$%\begin{equation}\label{mapr2}
r:\bnua{m}{2}(\alpha)\to\bbr^{m-1}-\{0\}\big/\{\pm 1\}\simeq\bbr P^{m-2}
$$%\end{equation}
This produces a cohomology class
$R\in H^1(\bnua{m}{2}(\alpha);\bbz_2)$ which
is the Stiefel-Whitney class of the double covering
$\nua{m}{2}(\alpha)\to\bnua{m}{2}(\alpha)$.
Lemma \ref{Ltotspace} holds true and Example \ref{excut<m>}
becomes:

\begin{Example}\label{excut<m>2}
\rm For the chamber $\langle m\rangle$, realized by
$a=(1,\dots,1,m-2)$,
the map $r$ is homotopic to a
diffeomorphism from $\nua{m}{2}(\langle m\rangle)$
onto $\bbr P^{m-3}$ and gives an identification of the
double covering
$\nua{m}{2}(\langle m\rangle)\to
\bnua{m}{2}(\langle m\rangle)$
with $S^{m-3}\to\bbr P^{m-3}$.
\end{Example}

If $\bar\cald(\alpha)$ denotes the total space of the $D^1$-bundle
associated to the double covering
$\nua{m}{2}(\alpha)\to\bnua{m}{2}(\alpha)$,
one proves, as in Proposition \ref{Pdouble}, that

\begin{Proposition}\label{Pdouble2}
If $\beta\in {\rm Ch}(\bbr_{\scriptscriptstyle\nearrow}^{m-1})$, then
$\bnua{m}{2}(\beta^+)$ is diffeomorphic to
the double of $\bar\cald(\beta)$.\cqfd
\end{Proposition}

\begin{Example}\label{excut<m1>}
\rm The chamber $\langle \{m,1\}\rangle=\langle m-1\rangle^+$  is
represented by\\ $a=(1/2,1,\dots,1,m-3)$.
As seen in Example \ref{excut<m>} one has that
$\nua{m}{3}(\langle m-1\rangle)$ is diffeomorphic to
$\bbc P^{m-4}$ and $\cald(\langle m-1\rangle)$ is
the disk bundle associated to the Hopf bundle.
Therefore $\nua{m}{3}(\langle \{m,1\}\rangle)$
is diffeomorphic to
$\bbc P^{m-3}\,\sharp\, \overline{\bbc P}^{m-3}$.
For planar polygons, one has
$\bnua{m}{2}(\langle m-1\rangle)$ is diffeomorphic to
$\bbr P^{m-4}$ and $\bar\cale(\langle m\rangle)
\to\bnua{m}{2}(\langle m\rangle)$ is
the double covering. Therefore,
$\bnua{m}{2}(\langle \{m,1\}\rangle)$
is diffeomorphic to
$\bbr P^{m-3}\,\sharp\,\overline{\bbr P}^{m-3}$
(of course, $\overline{\bbr P}^{m-3}=\bbr P^{m-3}$
when $m$ is even).
\end{Example}

\begin{ccote}\label{descrnongen}
Case where $\beta$ is non generic. \rm
When $\alpha=\beta^+$
with $\beta$ non generic, some partial information
about $\nua{m}{3}(\alpha)$ can still be gathered.
We proceed as in the proof of Proposition \ref{Pdouble},
with the same notations.
If $\varepsilon$ is small enough, the $(m-1)$-tuple
$b_{\delta}:=(b_2,\dots,b_m+\delta)$ is generic when
$0<|\delta|\leq\varepsilon$;
set $\beta_{\pm\varepsilon}:=\ch(b_{\pm\varepsilon})$.
The manifold $\check\cale(a)$ is now a cobordism between
between $\cale(\beta_{\varepsilon})$
and $\cale(\beta_{-\varepsilon})$.
By \cite[Thm 3.2]{Ha}, the map
$-\theta:\check\cale(a)\to [-\pi,0]$ is a Morse
function.
It has only one critical value,
the angle for which the diagonal length
$|\rho_m+\rho_1|$ is equal to $b_m$.
The preimage of this critical value
is diffeomorphic
to $\nua{m-1}{3}(\beta)$ and the (isolated) critical points
are the lined configurations.
There is one for each almost short subset
$I^=\in S_{m-1}(\beta)$ and its index
is equal to $2|I|$ (or $|I|$ for planar polygons).
As in Equation \eqref{lastdiffeo}, the space
$\nua{m}{3}(\alpha)$ is diffeomorphic to the quotient of
the cobordism $\check\cale(a)$ by the
following identifications on its two ends:
$\eta\sim g\cdot\eta$ for all $g\in SO(2)$,
when $\theta(\eta)=0,\pi$ (for planar polygons,
$g\in O(1)$).
\end{ccote}

As an application of the results of this section, we will
describe all the hexagon spaces in Section \ref{hexasp}.

\begin{ccote}\label{toric} Use of toric manifolds.\ \rm
Recall that a symplectic manifold $M^{2n}$ is called toric if it is endowed
with a Hamiltonian action of a torus $T$ of dimension $n$
(the maximal possible dimension for a Hamiltonian torus action).
The moment map $\mu:M\to {\rm Lie}(T)^*\approx\bbr^n$
has for image a convex polytope, the {\it moment polytope},
which determines
$M$ up to $T$-equivariant symplectomorphism (see \cite{Gu}).

The spatial polygon space $\nua{m}{3}(a)$
with its symplectic structure (see \ref{poisson})
may admit
Hamiltonian torus actions by so the so called {\it bending flows}
(see \cite{Kl}, \cite{KM}, \cite{HT}), which we recall now.
For $I\in\underline{m}$, define
$f_I:\nua{m}{3}(a)\to\bbr$ by
$f_I(\rho):=|\sum_{i\in I}\rho_i|$. If $f_I$ does not vanish,
it is a smooth map
which generates a Hamiltonian circle action
on $\nua{m}{3}(a)$. This action
rotates at constant speed
the set of vectors $\{\rho_i \mid i\in I\}$
around the axis $\sum_{i\in I}\rho_I$
(see \cite[\S\,2.1]{Kl}, \cite[Corollary 3.9]{KM}).
The non-vanishing of $f_I$ is equivalent to $I$ being
{\it lopsided}, that is there exists $i\in I$
with $a_i> \sum_{j\in I-\{i\}}a_j$ (see \cite{HT}).

Suppose that $\cali\subset\calp(\underline{m})$
is formed of lopsided subsets satisfying the following
``absorption condition'':
{\sl if $I,J\in\cali$ with $I\not=J$, then either
$I\cap J=\emptyset$ or one is contained in
the other}. Then, the Hamiltonian flow of the $f_i$'s of $I\in\cali$
commute and generate a Hamiltonian action of a torus $T_\cali$.
(see \cite[\S\,2.1]{Kl}, \cite[Lemma 2.1]{HT}).
Thus, when $\dim T_\cali=m-3$, the manifold
$\nua{m}{3}(a)$ is a toric manifold
which is determined by the moment polytope for the
$T_\cali$-action.

For example, when $m=5$, each chamber
$\alpha\in\ch(\bbr_{\scriptscriptstyle\nearrow}^5)$
has a representative
$a\in\alpha$ with $a_1\not=a_2$ and  $a_3\not=a_4$.
Therefore, $\nua{5}{3}(a)$ admits a Hamiltonian action
of the $2$-dimensional torus
$T_\cali$ for $\cali=\{\{1,2\},\{34\}\}$. This shows that
the diffeomorphism type of $\nua{5}{3}(a)$
is that of a toric manifold. The determination of all
the $2$-dimensional moment polytopes was the principle of the
classification of the $5$-gon spaces given in
\cite[\S\,6]{HK1}.

The same holds for $m=6$ since each chamber
$\alpha\in\ch(\bbr_{\scriptscriptstyle\nearrow}^6)$
has a representative
$a\in\alpha$ with $a_1\not=a_2$,
$a_3\not=a_4$ and $a_5\not=a_6$. Therefore, all
$\nua{6}{3}(\alpha)$ are diffeomorphic to
toric manifolds. The $3$-dimensional moment polytopes
can still be visualized but with more difficulties.

The above two cases generalizes in the following

\begin{Proposition}\label{protoric}
Let $\alpha\in\chrmcr$. Suppose that there exists
$a\in\alpha\cap\bbz^m$ with $a_m\geq\sum_{i=1}^{m-5}a_i$.
Then the diffeomorphism type of $\nua{m}{3}(\alpha)$
is that of a toric manifold.
\end{Proposition}

\preu One can find $a'\in\alpha$ arbitrarily close to $a$
so that $a_{m-4}\not=a_{m-3}$,
$a_{m-2}\not=a_{m-1}$ and $a_m>\sum_{i=1}^{m-5}a_i$.
Therefore the family of lopsided sets
$$\{m,1\}\, ,\, \{m,2,1\}\, ,\, \dots \, , \,
\{m,m-5,m-4,\dots,1\}\ ,\
\{m-4,m-3\}\, ,\, \{m-2,m-1\}$$
satisfy the absorption condition. Their bending flows
generate a Hamiltonian action of a torus of dimension
$m-3$, which proves the proposition.  \cqfd

\paragraph{Examples : } Consulting the table of the 135
seven-gons (see \cite{HRWeb}), we see that there are only three
$\alpha\in\ch(\bbr_{\scriptscriptstyle\nearrow}^7)$
for which $a=a_{\min}(\alpha)$ does not satisfy the hypothesis
of Proposition \ref{protoric}, that is, here, $a_7\geq a_1+a_2$.
These are
\begin{center}\small
\begin{tabular}{cc}
$\alpha$ & $a_{\min}(\alpha)$\\
\hline\rule{0pt}{3ex}
$\langle 754,762\rangle$ & (3,3,3,4,4,5,5)   \\[1mm]
$\langle 764\rangle$ &  (2,2,2,2,3,3,3)  \\[1mm]
$\langle 765\rangle$ &   (1,1,1,1,1,1,1)
\end{tabular}
\end{center}
\rm
Thus, all the other 133 heptagon spaces are diffeomorphic
to toric manifolds. We do not know whether the above
three heptagon spaces are diffeomorphic
to toric manifolds.

The same experiment with $m=8$ or $9$ gives the following
results: 217 elements of
$\rm{Ch}(\bbr_{\scriptscriptstyle\nearrow}^8)$ (out of 2400)
and 56550 elements of
$\rm{Ch}(\bbr_{\scriptscriptstyle\nearrow}^9)$ (out of 175428)
do not satisfy the hypothesis of \ref{protoric}.
\end{ccote}

\section{Cohomology invariants of $\nua{m}{3}(\alpha)$}\label{Scoh}

Let $\alpha$ be a chamber of $\rmpos$.
In \cite{HK2}, presentations of the cohomology
rings  $H^*(\nua{m}{3}(\alpha);\bbz)$
and  $H^*(\bnua{m}{2}(\alpha);\bbf_2)$
were obtained in terms of $\alpha$.
Our algorithms allowed us, with the help of a computer,
to find enough information about these rings
to prove that, for $5\leq m\leq 7$, $\alpha=\alpha'$ if and only if
the ${\rm mod}\, 2$ cohomology rings of $\nua{m}{3}(\alpha)$
and of  $\nua{m}{3}(\alpha')$ are isomorphic.

We start by the
Poincar\'e polynomial. Recall that $\nua{m}{3}(\alpha)$
has a cellular decomposition with only even-dimensional cells
\cite[\S\,4]{HK2}, so its Poincar\'e polynomial
is the same for any field $\bbf$:
$$P(t)=\sum_{i=0}^{2(m-3)} \dim_\bbf H^i(\nua{m}{3}(\alpha);\bbf)\, t^i$$
and has only terms of even degree. Moreover,
the polynomial $P(\sqrt{t})$
is the Poincar\'e polynomial of
$\bnua{m}{2}(\alpha)$ for the
coefficient field with two elements $\bbf_2$ \cite[\S\,9]{HK2}.
The first formula for computing $P(t)$ in terms of $\alpha$
was found by A. Kliachko \cite[Th. 2.2.4]{Kl}. We will use
the more economical formula, using only elements of
$S_m(\alpha)$, obtained in \cite[Cor. 4.3]{HK2}. With our notation,
this is:

\begin{Proposition}\label{Poincare}
Let $\alpha\in\chrmpos$.
The Poincar\'e polynomial of $\nua{m}{3}(\alpha)$ is
$$P(t)=
\frac{1}{1-t^2} \sum_
{J\in S_m(\alpha)}
(t^{2(|J|-1)}-t^{2(m-1-|J|)}).$$
%{J\in\calp(\underline{m-1})\atop J\cup\{m\}\in\alpha}{}
%(t^{2|J|}-t^{2(m-2-|J|)}).$$
\end{Proposition}

\begin{Remark}\label{chgnot}\rm
The difference between the formula in Proposition
\ref{Poincare} and that of \cite[Cor. 4.3]{HK2}
comes from that, there, the notation $S_m$ is used for the
set of $I\in \calp(\underline{m-1})$ such that
$J\cup\{m\}\in S$.
Recall that, here, $S_m=S\cap\calp_m(\underline{m})$.
So, each occurrence of $|J|$ in \cite[Cor. 4.3]{HK2}
is replaced here by $|J|-1$.
\iffalse
The indexing set for the sum in Proposition \ref{Poincare}
is the same as in \cite[Cor. 4.3]{HK2} but expressed
differently.
In \cite{HK2}, the notation $S_m$ is used for this set,
%set of $I\in \calp(\underline{m-1})$ such that
%$J\cup\{m\}\in S$.
which does not match with our notation.
Recall that, here, $S_m=S\cap\calp_m(\underline{m})$.
\fi
\end{Remark}

For $\beta\subset\calp(\underline{m})$,
denote by $\ns_i(\beta)$ the number of sets $I\in\beta$
with $|I|=i+1$. The formula of Proposition \ref{Poincare}
gives the Betti numbers
$b_{2i}:=b_{2i}(\alpha):=\dim_\bbf H^i(\nua{m}{3}(\alpha);\bbf)$
as the solution of the system of equations
\stepcounter{dummy}\begin{equation}\label{eqbi}
b_{2i}-b_{2i-2}= \ns_i(S_m(\alpha))-\ns_{m-2-i}(S_m(\alpha)),
\end{equation}
starting with $b_{2i}=0$ if $i<0$. For instance, if
$\alpha=\langle 54\rangle$, realized by $(1,1,1,1,1)$, one has
$$S_5(\alpha)=\{5,51,52,53,54\}$$
thus $\ns_0(S_5(\alpha))=1$, $\ns_1(S_5(\alpha))=4$ and the other
$\ns_i(S_5(\alpha))$ vanish.
This gives $b_0=1$, $b_2=5$ and $b_4=1$, which are
indeed the Betti numbers of $\nua{5}{3}(\langle 54\rangle)
=\bbc P^2\sharp 4 \overline{\bbc P^2}$ \cite[Example 10.4]{HK2}.

As our computer algorithm had to list all the sets of $S_m(\alpha)$
(for instance, for the realization), the numbers $\ns_i(S_m(\alpha))$
are available and so are the $b_{2i}$'s.

We now recall the presentation of $H^*(\nua{m}{3}(\alpha);\bbz)$
obtained in \cite[Thm 6.4]{HK2}. Taking care of Remark \ref{chgnot},
this gives:

\begin{Proposition}\label{prescoh}
Let $\alpha\in\chrmpos$.
The cohomology ring of the polygon space
$\nua{m}{3}(\alpha)$ with coefficient in a ring $\Lambda$
 is
$$
H^*(\nua{m}{3}(\alpha);\Lambda)=
\Lambda[R,V_1,\dots,V_{m-1}]/\cali(\alpha)$$
where $R$ and $V_i$ are of degree 2
and $\cali(\alpha)$ is the ideal
of $\Lambda[R,V_1,\dots,V_{m-1}]$
generated by the three families\smallskip

\begin{tabbing} \renewcommand{\arraystretch}{5}
\kern .7 truecm \= \+ (R1)\quad   \=
$R^2\sum_{S\subset L \atop S\in S_m(\alpha)} \big(\prod_{i\in S}V_i\big)
R^{|L-S|-1}$  \=    \= \kill

(R1)\> $V_i^2+ RV_i$ \>  \> $i=1,\dots ,m-1$  \\ \bigskip

(R2)  \> ${\displaystyle \prod^{\ }_{i\in L} V_i}$  \>\>
 for all $L\in\calp(\underline{m-1})$ with $L\cup\{m\}
\notin S(\alpha)$
  \\ \bigskip

 (R3)  \> ${\displaystyle \sum^{\ }_{S\subset L \atop
S\cup\{m\}\in\alpha}
\big(\prod_{i\in S}V_i\big) R^{|L-S|-1}}$
\>\> for all $L\in\calp(\underline{m-1})$ with $L\notin S(\alpha)$
\end{tabbing}
\end{Proposition}

We shall use this presentation, first to compute a homotopy invariant
$r_\cup(\alpha)\in\bbn$, defined as the rank of the linear map
$x\mapsto x\cup x$ from
$H^2(\nua{m}{3}(\alpha);\bbf_2)$ to
$H^4(\nua{m}{3}(\alpha);\bbf_2)$
(recall that $x\mapsto x^2$ is a linear map in an
algebra over $\bbf_2$).
Observe that $r_\cup(\alpha)$ is also the rank
of the same map from
$H^1(\bnua{m}{2}(\alpha);\bbf_2)$ to
$H^2(\bnua{m}{2}(\alpha);\bbf_2)$.

\begin{Proposition}\label{rcup}
For all $\alpha$ a chamber of $\rmcr$, one has
\begin{equation}\label{rcupfor}
r_\cup(\alpha)=1+\ns_1(S_m(\alpha))-\ns_{m-3}(S_m(\alpha))-\ns_{m-4}(S_m(\alpha)).
\end{equation}
\end{Proposition}

\preu
One has $1+\ns_1(S_m(\alpha))-\ns_{m-3}(S_m(\alpha))=b_2$
by equations \ref{eqbi}. Let us first consider the case
$\ns_{m-3}(S_m(\alpha))\not=0$. This means that $S_m(\alpha)$
contains a set with $m-2$ elements and thus contains
the smallest of those
for the order $\fpp$, which is
$I:=\{m,m-3,m-4,\dots,1\}$.
Then $\alpha=\langle I\rangle$.
Indeed, if $\alpha\not=\langle I\rangle$, then
$\alpha$ would contain $J:=\{m,m-2\}$, which is impossible
since $\bar{J}\fpp I$. Therefore,
$$\ns_{1}(S_m(\alpha))=m-3 \quad ,\quad
\ns_{m-4}(S_m(\alpha))=m-3  \quad ,\quad
\ns_{m-3}(S_m(\alpha))=1$$
and the right hand member of
\eqref{rcupfor} is equal to zero. On the other hand, Relator (R3)
of Proposition \ref{prescoh}, with $L=\{m-2,m-1\}$ and $S=\emptyset$
gives the equality $R=0$. By Relators (R1), all squares vanish
and $r_\cup=0$. Formula \eqref{rcupfor} is then proven in the
case $\ns_{m-3}(S_m(\alpha))\not=0$. Observe that
$\nua{m}{3}(\langle I \rangle)$
is diffeomorphic to a product of $m-3$ copies of the sphere $S^2$
\cite[Example 10.2]{HK2}.

Assume then that $\ns_{m-3}(S_m(\alpha))=0$. By Proposition \ref{prescoh}, the vector space\\
\kern 15mm
$H^2(\nua{m}{3}(\alpha);\bbf_2)
=H^2(\nua{m}{3}(\alpha);\bbz)\otimes\bbf_2$
has the basis $R,V_1,\dots,V_p$ for $p=\ns_1(S_m(\alpha))$.
Indeed, $V_i=0$ for $i>p$ by Relator (R2) with $L=\{i\}$.
The image of $x\mapsto x^2$ is generated by
$R^2,V_1^2,\dots,V_p^2$. The relations between these generators
come from relators (R3) of Proposition \ref{prescoh}
with $|L|=3$. For such an $L=\{i,j,k\}$ (denoted by $ijk$), the relation is
\stepcounter{dummy}\begin{equation}\label{relquad}
R^2+RV_i+RV_j+RV_k + V_iV_j + V_iV_k + V_jV_k = 0 .
\end{equation}
The three last terms of the left hand member
vanish by Relator (R2). Indeed, since
$ijk\notin\alpha$  and $ijk\fpp ijm$, then $ijm\notin\alpha$.
By Relator (R1), Equation \eqref{relquad} becomes
\stepcounter{dummy}\begin{equation}\label{eqijk}
R^2-V_i^2-V_j^2-V_k^2=0.
\end{equation}
To establish Proposition \ref{rcup}, it is enough to prove that,
for distinct $i,j,k$, Equations \eqref{eqijk}
are independent. But, by Proposition \ref{prescoh},
The vector space $H^4(\nua{m}{3}(\alpha);\bbf_2)$
is generated by the $b_1$ elements
$R^2,V_1^2,\dots,V_p^2$ together with the
$\ns_2(S_m(\alpha))$ non-vanishing
products $V_iV_j$ ($ijm\in\alpha$) and these generators
are just subject to Equations \eqref{eqijk}. This implies that
$b_2\geq b_1 + \ns_2(S_m(\alpha)) - \ns_{m-4}(S_m(\alpha))$.
By Equations \eqref{eqbi},
this inequality is an equality, showing that Equations \eqref{eqijk}
are independent.  \cqfd

The idea of our last cohomology invariant
Was given to us by R. Bacher.
Define $\cali_k(\alpha)$ to be the
ideal of $\bbf[R,V_1,\dots,V_{m-1}]$ generated by the elements
of $\cali(\alpha)$ which are polynomials of degree $\leq k$
in the variables $R$ and $V_i$'s (then giving elements
of degree $\leq 2k$ in $H^*(\nua{m}{3}(\alpha);\bbf)$).
Define $\sol_k(\alpha;\bbf)\subset \bbf P^{m-1}$
to be the projective variety defined by
the equations $W=0$ for all $W\in\cali_2(\alpha)$.

\begin{Proposition}\label{projvar}
a) Let $\alpha,\alpha'$ be chambers of $\rmcr$.
Any graded ring isomorphism from
$H^*(\nua{m}{3}(\alpha);\bbf)$ onto
$H^*(\nua{m}{3}(\alpha');\bbf)$ induces, for $k\geq 1$,
a bijection from
$\sol_k(\alpha';\bbf)$ to $\sol_k(\alpha;\bbf)$.

b) Suppose that $k\geq 2$.
Then, any element $\zeta\in\sol_k(\alpha;\bbf)$
has a unique representative
of the form $(-1,v_1,\dots,v_{m-1})\in\bbf^m$ with $v_i=0$ or $1$.
In particular, the set $\sol_k(\alpha;\bbf)$ is finite.

c) The finite set $\sol_2(\alpha;\bbf)$ does not depend on the field $\bbf$.
\end{Proposition}

\preu

\noindent\it Proof of a) : \rm
Let $q:H^*(\nua{m}{3}(\alpha);\bbf)
\to H^*(\nua{m}{3}(\alpha');\bbf)$
be a graded ring homomorphism. By Proposition \ref{prescoh},
the homomorphism $q$ is covered by a graded ring
homomorphism
$\tilde q:\bbf [R,V_1,\dots,V_{m-1}]\to\bbf [R',V_1',\dots,V_{m-1}']$
which sends $\cali_k(\alpha)$ into $\cali_k(\alpha')$ for all $k$.
Such a lifting $\tilde q$ is well defined up to a homomorphism
with image in $\cali_1(\alpha')$, therefore $q$
functorialy induces homomorphisms
$\tilde q:\bbf [R,V_1,\dots,V_{m-1}]/\cali_k(\alpha)
\to\bbf [R',V_1',\dots,V_{m-1}']/\cali_k(\alpha')$ for all $k\geq 1$.
Observe that $\sol_k(\alpha;\bbf)$ can be identified
with the projectivization of the vector space of ring homomorphisms
from  $\bbf [R,V_1,\dots,V_{m-1}]/\cali_k(\alpha)$  to $\bbf$.
Therefore, the homomorphism $q$ will
functorialy induce maps
$\hat q:\sol_k(\alpha';\bbf)\to \sol_k(\alpha;\bbf)$
for all $k\geq 1$. By this functoriality, if $q$ is an isomorphism,
then $\hat q$ is a bijection, which proves a).

\noindent\it Proof of b) : \rm
let $z=(r,v_1,\dots,v_{m-1})\in\bbf^m-\{0\}$ represent an
element $\zeta\in \sol_k(\alpha;\bbf)$. Then,  $r\not=0$,
since, otherwise relators
(R1) of Proposition \ref{prescoh} would give
equations $v_i^2=-rv_i$,
implying that $z=0$.
Then, $[z]$ has a unique representative with $r=-1$
and the equations  $v_i^2=-rv_i$
imply that $v_i\in\{0,1\}$.

\noindent\it Proof of c) : \rm
The equations defining $\sol_2(\alpha;\bbf)$,
coming from Relators (R1)-(R3) of \ref{prescoh}, are,
for $i,j,k=1\dots,m-1$:
\begin{enumerate}
\renewcommand{\labelenumi}{(\roman{enumi})}
\item $v_i^2=-rv_i$. Having normalized $r=-1$,
these are equivalent to $v_i=0,1$.
\item $v_i=0$, for $\{i,m\}\notin\alpha$ and $v_iv_j=0$
for  $\{i,j,m\}\notin\alpha$.
\item Equations \eqref{eqijk} which, after (i), become
$v_i^2+v_j^2+v_k^2=1$, for $\{i,j,k\}\notin \alpha$.
\end{enumerate}
The solutions $v_i=0,1$ of Equations (ii) are clearly independent
of the ground field $\bbf$. The solutions $v_i=0,1$ of
an equation like $v_i^2+v_j^2+v_k^2=1$
seem, a priori, to depend on the characteristic of $\bbf$. But,
as seen just before Equation \eqref{eqijk}, such an equation
occurs only if $r_ir_j=r_jr_k=r_kr_i=0$. Thus, Equation (iii)
is equivalent to the fact that exactly one of the $v_i,v_j,v_k$
is equal to one, a condition independent of $\bbf$.
\cqfd

\begin{Definition}
We set $s(\alpha)$ to be the number of elements of\
$\ \sol_2(\alpha;\bbf)$.
This does not depend on the field $\bbf$ by Proposition \ref{projvar}.
\end{Definition}

It is not difficult to compute $s(\alpha)$ with or without
the help of a computer. We select
the elements of $W(R,V_1,\dots,V_{m-1})\in \cali_2(\alpha)$
and count how many of them vanish when $R=-1$ and $V_i\in\{0,1\}$.

\begin{ccote}\label{pfprB} Proof of Proposition B: \rm
The following implies Proposition B of the introduction:

\begin{Proposition}\label{chmod2disting}
Let $5\leq m\leq7$ and let $\alpha,\alpha'\in\chrmcr$. Then
$\alpha=\alpha'$ if and only if $\nua{m}{3}(\alpha)$
and $\nua{m}{3}(\alpha')$ have the same Betti numbers,
$r_\cup(\alpha)=r_\cup(\alpha')$ and $s(\alpha)=s(\alpha)'$.
\end{Proposition}

\preu  For $m=5$, by the list of Table V,
The only case where two $5$-gon spaces have the same
Betti numbers are
$\nua{5}{3}(\langle 52\rangle)
\approx \bbc P^2\,\sharp\,\overline{\bbc P}^2$
and  $\nua{5}{3}(\langle 521\rangle)\approx
S^2\times S^2$. But $r_\cup(\langle 52\rangle)=1$ while
$r_\cup(\langle 521\rangle)=0$. (Taking $\bbf_2$ coefficients
is important here: for instance,
these two spaces have isomorphic
cohomology ring with real coefficients.)

For $m=6$, the list of Table VI in Section \ref{hexasp}
has been sorted by lexicographic order of the triple
$b_2(\alpha),r_\cup(\alpha),s(\alpha)$.
One thus can check that no such triples occur twice.
The same holds for $m=7$ with
$b_2(\alpha),b_4(\alpha),r_\cup(\alpha),s(\alpha)$
(table in \cite{HRWeb}).
\end{ccote}

\begin{Remark}\rm By \cite[\S\,9]{HK2},
the cohomology ring $H^*(\bnua{m}{2}(\alpha);\bbf_2)$
admits the presentation
of Proposition \ref{prescoh}, with $R$ and the $V_i$'s
of degree $1$. Therefore, the above invariants are
${\rm mod\,} 2$ cohomology invariants of the spaces
$\bnua{m}{2}(\alpha);\bbf_2$  and Proposition
\ref{chmod2disting} holds true for these spaces.
\end{Remark}

\begin{Problem}\label{cohnonrea}\rm
When a virtual genetic code $\gamma\in\calg_m$
is not realizable (for instance when $m=9$), it
gives rise as well
to a non-trivial graded ring.
Is this ring the cohomology ring of a space?
Does it satisfy Poincar\'e duality?
Is this ring the cohomology ring of a manifold?
\end{Problem}

\section{The hexagon spaces}\label{hexasp}

As for Tables III-V of Section \ref{Scalna}
the first column of Table VI below
contains the list of the 21 genetic codes of
type $6$, all realized by a chamber $\alpha$
whose minimal realization
$a_{\rm min}(\alpha)$, using conventional
representatives (see \ref{convzeros}),
is written in the second column.
The next three columns give the cohomology invariants
of $\nua{6}{3}(\alpha)$ or $\bnua{6}{2}(\alpha)$
which are defined in Section \ref{Scoh}, with the
abbreviations
$$
b:=\dim_\bbf H^2(\nua{6}{3}(\alpha);\bbf)
= \dim_{\bbf_2} H^1(\bnua{6}{2}(\alpha);\bbf_2)\ , \
r_{\cup} =r_{\cup}(\alpha) \ ,\ s=s(\alpha).$$
By Poincar\'e duality, the number $b$ determines
the Poincar\'e polynomial $P(t)$ of
$\nua{6}{3}(\alpha)$ which is
$$P(t)=1+bt^2+bt^4+t^6$$
(for $\bnua{6}{2}(\alpha)$, this would be
$P(t)=1+bt^+bt^2t^3$).
The $6$-gon
spaces have been listed by the lexicographic order of
the triples $(b,r_\cup,s)$, showing that
the homotopy type of the hexagon spaces
in $\bbr^2$ or $\bbr^3$ are distinguished by
these cohomology invariants.

The last two columns contain some geometric descriptions
of $\bnua{6}{2}(\alpha)$ and $\nua{6}{2}(\alpha)$
obtained by the methods discussed in Section \ref{Scalna}.
In the last column, we see the hexagon spaces coming
from the 7 generic pentagons by adding a tiny vector; their
descriptions uses Proposition \ref{timesS1}.
Lines 2 and 3, illustrate Examples \ref{excut<m>2}
and \ref{excut<m1>}. Line 4 is a special case of
\cite[Example 10.2]{HK2} (see also the proof of
Proposition \ref{rcup}). The other descriptions come
from the method of \ref{descrnongen}. Using the notations
of \ref{descrnongen}, the $3$-manifold
$\check\cale(a)$ is a cobordism between the orientable
surfaces $\cale(\beta_{\varepsilon})$
and $\cale(\beta_{-\varepsilon})$. When they are connected
(all cases except Line 7), we denote their genus
respectively by $g_+$ and $g_-$. The map
$-\theta:\check\cale(a)\to [-\pi,0]$ is a Morse
function with $n_1$ critical points of index $1$
and $n_2$ critical points of index $2$. This situation
is indicated in Table VI by the writing $[g_+;n_1,n_2;g_-]$
(observe that $g_-=g_+ +n_1-n_2$).
The orientable $3$-manifold
$\bnua{6}{2}(\alpha)$ is diffeomorphic to the quotient of
the cobordism $\check\cale(a)$ by, on each end,
a free involution reversing the orientation.

In other words, the notation $[g_+;n_1,n_2;g_-]$
tells us that the orientable $3$-manifold
$\bnua{6}{2}$ is obtained in the following way.
For $\Sigma$ a surface, denote by $\bar\cald(\Sigma)$ the mapping
cylinder of the orientation covering of $\Sigma$.
Let $W_+=\bar\cald(\su{g_+})$ union with $n_1$ $1$-handles
and $W_+=\bar\cald(\su{g_-})$ union with $n_2$ $1$-handles.
If we require that $W_\pm$ are orientable, they
are well defined since there is
only one way, up to diffeomorphism isotopic to the
identity, to attach $1$-handles to $\Sigma^{or}_{g_\pm}\times [0,1]$
in order to obtain an orientable manifold.
Thus, $\bnua{6}{2}$ is obtained by gluing $W_+$ to
$W_-$ by a diffeomorphism of their boundary. In the case
where $n_1=n_2=0$, one has $W_+=W_-$ and
Proposition \ref{Pdouble2} says that
the gluing diffeomorphism is the identity.
We were not able to identify this gluing diffeomorphism
in the other cases, so, {\it a priori}, the numbers
$[g_+;n_1,n_2;g_-]$ do not determine the homeomorphism
type of $\bnua{6}{2}$.

In the case $\alpha=\langle 632\rangle$ (line 7 of the table),
$\bnua{5}{2}(\beta_\varepsilon)=T^2$
(the only case where it is orientable). Therefore
$W_+=T^2\times [-1,1]$ and $W_-=\bar\cald(\su{2})$.

\vskip 5mm
{\small
\begin{center}
\setcounter{compol}{0}
\begin{supertabular}{rccccccccc}
\multicolumn{8}{c}{\bf Table VI : the $6$-gon spaces}\\
\hline \rule{0pt}{3ex}
& $\alpha$ & $a_{\rm min}(\alpha)$ & $b$ & ${\rm r}_\cup$ & $s$ &
%\hspace{1pt}$\nua{6}{3}(\alpha)$  &
\hspace{1pt} $\bnua{6}{2}(\alpha)$ \kern 5mm &
$\nua{6}{2}(\alpha)$
\\[3pt]\hline
\rule{0pt}{3ex}
\pvi{0}{0}{0}{0,0,0,0,0,1}{1}{1}{}\ccvi{$\emptyset$}{}{$\emptyset$}{$\emptyset$}%x
\pvi{1}{1}{1}{1,1,1,1,1,4}{9}{2}{6}\ccvi{$\bbc P^3$}{}{$\bbr P^3$}{$S^3$}%x
\pvi{2}{2}{2}{0,1,1,1,1,3}{7}{3}{61}%
\ccvi{$\bbc P^3\,\sharp\,\overline{\bbc P}^3$}{}%
{$\bbr P^3\,\sharp\,\overline{\bbr P}^3$}{$S^2\times S^1$}%x
\pvi{3}{0}{0}{0,0,0,1,1,1}{3}{9}{6321}%
\ccvi{$(S^2)^3$}{}{$T^3$}{$T^3{\scriptstyle\coprod}T^3$}%x
\pvi{3}{2}{0}{0,0,1,1,1,2}{5}{5}{621}
\ccvi{}{}{$[1;0,0;1]$}{$T^3$}%x
\pvi{3}{3}{3}{1,1,2,2,2,5}{13}{4}{62}
\ccvi{}{}{$[0;1,0;1]$}{}%x
\pvi{4}{1}{1}{1,1,1,3,3,4}{13}{8}{632}
\ccvi{}{}{{\tiny see above}}{}%x
\pvi{4}{2}{0}{0,1,1,2,2,3}{9}{7}{631}
\ccvi{}{}{$[2;0,0;2]$}{$\Sigma_2^{or}\times S^1$}%x
\pvi{4}{3}{1}{1,1,2,3,3,5}{15}{14}{621,63}
\ccvi{}{}{$[1;1,0;2]$}{}%x
\pvi{4}{4}{4}{1,1,1,2,2,4}{11}{6}{63}
\ccvi{}{}{$[0;2,1;2]$}{}%x
\pvi{5}{2}{0}{0,1,1,1,2,2}{7}{11}{641}
\ccvi{}{}{$[3;0,0;3]$}{$\Sigma_3^{or}\times S^1$}%x
\pvi{5}{2}{2}{1,1,1,2,3,3}{11}{19}{632,64}
\ccvi{}{}{$[2;1,1;2]$}{}%x
\pvi{5}{3}{1}{1,2,2,3,4,5}{17}{17}{631,64}
\ccvi{}{}{$[2;1,0;3]$}{}%x
\pvi{5}{4}{2}{1,1,2,2,3,4}{13}{15}{621,64}
\ccvi{}{}{$[1;2,0;3]$}{}%x
\pvi{5}{5}{5}{1,1,1,1,2,3}{9}{10}{64}
\ccvi{}{}{$[0;3,0;3]$}{}%x
\pvi{6}{2}{0}{0,1,1,1,1,1}{5}{13}{651}
\ccvi{}{}{$[4;0,0;4]$}{$\Sigma_4^{or}\times S^1$}%x
\pvi{6}{3}{1}{1,2,2,2,3,3}{13}{21}{641,65}
\ccvi{}{}{$[3;1,0;4]$}{}%x
\pvi{6}{3}{3}{1,1,1,2,2,2}{9}{20}{632,65}
\ccvi{}{}{$[2;2,1;3]$}{}%x
\pvi{6}{4}{2}{1,2,2,3,3,4}{15}{18}{631,65}
\ccvi{}{}{$[2;2,0;4]$}{}%x
\pvi{6}{5}{3}{1,1,2,2,2,3}{11}{16}{621,65}
\ccvi{}{}{$[1;3,0;4]$}{}%x
\pvi{6}{6}{6}{1,1,1,1,1,2}{7}{12}{65}
\ccvi{}{}{$[0;4,0;4]$}{}%x
\hline
\end{supertabular}
\end{center}
}\rm

\begin{ccote} Problem. \rm  Describe more precisely the $3$-dimensional
manifolds $\bnua{6}{2}$, for instance in terms of the Kirby calculus.
\end{ccote}

%\sk{3}\goodbreak

\sk{0}
{\footnotesize
Section de Math\'ematiques\hfill hausmann{@}math.unige.ch\\
Universit\'e de Gen\`eve, B.P. 240 \hfill rodriguez{@}math.unige.ch\\
CH-1211 Geneva 24,
Switzerland
}

\end{document}